\documentstyle[10pt]{amsart}

\newtheorem{prop}{Proposition}

\newtheorem{th}{Theorem}
\newtheorem{lemma}{Lemma}
\newtheorem{cor}{Corollary}

\title[Compatible structures  and their Gray-Hervella classes \hfill]{Compatible almost complex structures on twistor spaces and their Gray-Hervella classes}
\author{Danish Ali, Johann Davidov and Oleg Mushkarov}
\address{Danish Ali\\Abdus Salam School of Mathematical Sciences\\
GC University, Lahore, Pakistan}
\email{danish$\_$mathali@@yahoo.com}
\address{Johann Davidov\\ Institute of Mathematics and Informatics \\
Bulgarian Academy of Sciences\\ Acad. G.Bonchev Str. Bl.8\\ 1113
Sofia\\ Bulgaria and \\ Abdus Salam School of Mathematical Sciences\\
GC University, Lahore, Pakistan}\email{jtd@@math.bas.bg}

\address{Oleg Mushkarov \\Institute of Mathematics and Informatics \\
Bulgarian Academy of Sciences\\ Acad. G.Bonchev Str. Bl.8\\ 1113
Sofia\\ Bulgaria and \\Abdus Salam School of Mathematical Sciences\\
GC University, Lahore, Pakistan} \email{muskarov@@math.bas.bg}

\begin{document}

\begin{abstract}
In this paper we determine the Gray-Hervella classes of the
compatible almost complex structures on the twistor spaces of
oriented Riemannian four-manifolds considered  by G. Deschamps in
\cite {De}.
\medskip

\vspace{0,1cm}
\noindent
2000 {\it Mathematics Subject Classification}. Primary 53C15, 53C25.

\vspace{0,1cm} \noindent {\it Key words: twistor spaces, compatible
almost complex structures, Gray-Hervella classes,}
\end{abstract}

\thispagestyle{empty}

\maketitle
\vspace{0.5cm}

\section{Introduction}

In the 80's of the last century A. Gray and L.M. Hervella \cite
{GH} have proposed a natural classification of almost Hermitian
manifolds by studying a representation of the unitary group on the
space of tensors satisfying the same identities as the covariant
derivative of the K\"ahler form of an almost Hermitian manifold.
This representation has four irreducible components, which
determine sixteen classes of almost Hermitian manifolds playing an
important role in Hermitian geometry.

The main purpose of the present paper is to determine the
Gray-Hervella classes of the so-called compatible almost Hermitian
structures on the twistor space of an oriented four-dimensional
Riemannian manifold $(M,g)$ considered by G. Deschamps in \cite
{De}. The (positive) twistor space of $(M,g)$ is the total space of
the $2$-sphere bundle $\pi:{\cal Z}\to M$ consisting of all unit
$(+1)$-eigenvectors of the Hodge star operator acting on
$\Lambda^2TM$. The fibre of ${\cal Z}$ at a point $p\in M$ can be
identified with space of all complex structures on the tangent space
$T_pM$ compatible with the metric and orientation of $M$. The
Levi-Civita connection of $(M,g)$ gives rise to a splitting $T{\cal
Z}={\cal H}\oplus{\cal V}$ of the tangent bundle of ${\cal Z}$ into
horizontal and vertical subbundles. Then, following
Atiyah-Hitchin-Singer \cite{AHS} and Eells-Salamon \cite{ES}, one
can define two almost complex structures ${\cal J}_1$ and ${\cal
J}_2$ on the six-manifold ${\cal Z}$ as follows. On every horizontal
space ${\cal H}_{\sigma}$, $\sigma\in{\cal Z}$, ${\cal J}_1$ and
${\cal J}_2$ are both the horizontal lift of the complex structure
$\sigma:T_{\pi(\sigma)}M\to T_{\pi(\sigma)}M$. The vertical space
${\cal V}_{\sigma}$ is the tangent space at the point $\sigma$ of
the fibre of ${\cal Z}$ through $\sigma$. This fibre is a unit
$2$-sphere and  ${\cal J}_1$  is defined on ${\cal V}_{\sigma}$ as
the standard complex structure of the sphere, while ${\cal J}_2$ is
the conjugate complex structure $-{\cal J}_1$. By a famous theorem
of Atiyah-Hitchin-Singer \cite{AHS} the almost complex structure
${\cal J}_1$ is integrable (i.e. comes from a complex structure) if
and only if the base manifold $(M,g)$ is anti-self-dual. In
contrast, the almost complex structure ${\cal J}_2$ is never
integrable by a result of Eells-Salamon \cite{ES} but it is very
useful in harmonic maps theory. The almost complex structures ${\cal
J}_1$ and ${\cal J}_2$ are compatible with the $1$-parameter family
of Riemannian metrics $h_t=\pi^{\ast}g+tg^{v}$, $t>0$, where $g^{v}$
is the restriction to the fibres of ${\cal Z}$ of the metric on
$\Lambda^2TM$ induced by $g$. The Gray-Hervella classes of the
almost Hermitian structures $(h_t,{\cal J}_1)$ and $(h_t,{\cal
J}_2)$ have been determined in \cite{M}.

It has been  observed by G. Deschamps in \cite{De} that one can
obtain other almost complex structures on ${\cal Z}$ compatible with
the metrics $h_t$ by means of a fibre-preserving map $f:{\cal
Z}\to{\cal Z}$. Given such a map, the corresponding almost complex
structure ${\cal J}_f$ on ${\cal H}_{\sigma}$ is defined as the
horizontal lift of $f(\sigma)$; on ${\cal V}_{\sigma}$ it is set to
be equal to ${\cal J}_1$. Thus, if $f=id$ we obtain the almost
complex structure ${\cal J}_1$ and if $f$ is the antipodal map
$\sigma\to -\sigma$ we get $-{\cal J}_2$. In Section 3 of the
present paper we derive coordinate-free formulas for the covariant
derivative of the K\"ahler $2$-form of the almost Hermitian
structure $(h_t,{\cal J}_f)$ corresponding to an arbitrary $f$. We
use these formulas to determine the Gray-Hervella classes of
$(h_t,{\cal J}_f)$ for some particular fibre-preserving maps
$f:{\cal Z}\to{\cal Z}$. More precisely, let $(M,g,J)$ be an almost
Hermitian manifold of real dimension four. Consider $M$ with the
orientation induced by the almost complex structure $J$. Then $J$ is
a section of the (positive) twistor space ${\cal Z}$ of $(M,g)$. At
any point $p\in M$, the complex structure $J_p:T_pM\to T_pM$ is a
point of the fibre ${\cal Z}_p$ at $p$ of ${\cal Z}$. Take a complex
number $\lambda$. Since ${\cal Z}_p$ is the unit $2$-sphere, we can
compose the stereographic projection of $(h_t,{\cal J}_f)$ from the
point $J_p$ with the linear map $z\to\lambda z$ of the complex
plane, then go back to the sphere by the inverse stereographic
projection. In this way we obtain a fibre-preserving map
$f^{+}_{\lambda}:{\cal Z}\to{\cal Z}$ whose restriction to every
fibre is a holomorphic map. If we use in a similar way the
stereographic projection from the point $-J_p$, we get a map
$f^{-}_{\lambda}$ whose restrictions to the fibres of ${\cal Z}$ are
anti-holomorphic maps. In particular, $f_1^+(\sigma)=\sigma$, so
${\cal J}_{f_1^+}={\cal J}_1$, the Atiyah-Hitchin-Singer almost
complex structure, whereas $f_{-1}^-(\sigma)=-\sigma$ and ${\cal
J}_{f_{-1}^-}=-{\cal J}_2$, the conjugate structure of the
Eells-Salamon almost complex structure. For $\lambda=0$, we have
$f_0^\pm\equiv\mp J$; note that the structures $J$ and $-J$ induce
the same orientation since $dim\,M=4$ and belong to the same
Gray-Hervella classes. In the case $\lambda=0$, the integrability
condition for the corresponding almost complex structure ${\cal
J}_{f_0^\pm}$ has been given in \cite{De}(where this structure is
denoted by ${\Bbb J}_{\infty}$). In Section 4, Theorem~\ref{GH} we
establish all possible Gray-Hervella classes of the almost Hermitian
structure $(h_t,{\cal J}_{f_0^\pm})$ on the twistor space ${\cal Z}$
and found the geometric conditions on the base manifold $M$ under
which this structure belongs to each of these classes. In the case
when $\lambda\neq 0,1$ and the base manifold $M$ is K\"ahler, the
integrability condition for the almost complex structure ${\cal
J}_{f^{+}_{\lambda}}$ has been obtained in \cite{De} (the structure
being denoted there by ${\Bbb J}_{\lambda Id}$ ). Under the
assumptions that $\lambda\neq 0$ and $M$ is K\"ahler, in Section 5,
Theorem~\ref{GHC1} we determine the Gray-Hervella classes of the
almost Hermitian structure $(h_t,{\cal J}_{f^{\pm}_{\lambda}})$. At
the end of this section, we discuss also the case when $|\lambda|=1$
without the K\"ahler assumption on the base manifold $M$. In order
to keep the length of the paper reasonable, we discuss only some of
the basic Gray-Hervella classes.

\section{Preliminaries}

Let $(M,g)$ be an oriented  Riemannian manifold of dimension four.
The metric $g$ induces a metric on the bundle of two-vectors
$\pi:\Lambda^2TM\to M$ by the formula
$$
g(v_1\wedge v_2,v_3\wedge v_4)=\frac{1}{2}det[g(v_i,v_j)].
$$
The Levi-Civita connection of $(M,g)$ determines a connection on the
bundle $\Lambda^2TM$, both denoted by $\nabla$, and the
corresponding curvatures are related by
$$
R(X\wedge Y)(Z\wedge T)=R(X,Y)Z\wedge T+Z\wedge R(X,Y)T
$$
for $X,Y,Z,T\in TM$. Let us note that we adopt the following
definition for the curvature tensor $R$ :
$R(X,Y)=\nabla_{[X,Y]}-[\nabla_{X},\nabla_{Y}]$. Then the curvature
operator ${\cal R}$ is the self-adjoint endomorphism of
$\Lambda^2TM$ defined by
$$
g({\cal R}(X\wedge Y),Z\wedge T)=g(R(X,Y)Z,T).
$$

\ The Hodge star operator defines an endomorphism $\ast$ of
$\Lambda^2TM$ with $\ast^2=Id$. Hence we have the decomposition
$$
\Lambda^2TM=\Lambda^2_{-}TM\oplus\Lambda^2_{+}TM
$$
where $\Lambda^2_{\pm}TM$ are the subbundles of $\Lambda^2TM$
corresponding to the $(\pm 1)$-eigenvalues of the operator $\ast$.

For every $p\in M$, the group $SO(4)$ acts in a natural way on the
space of $4$-tensors on $T_pM$ having the same symmetries as the
Riemannian curvature tensor. The irreducible decomposition of this
space under the action of $SO(4)$, found by Singer and Thorpe
\cite{ST}, gives the following decomposition of the curvature
operator
\begin{equation}\label{cur-dec}
{\cal R}=\frac{s}{6}Id+{\cal B}+{\cal W}_{+}+{\cal W}_{-}
\end{equation}
where $s$ is the scalar curvature, the operator ${\cal B}$
represents the traceless Ricci tensor and ${\cal W}_{\pm}$ are the
restrictions on $\Lambda_{\pm}TM$ of the operator ${\cal W}$
corresponding the the Weyl conformal tensor. These operators are
symmetric and  ${\cal B}$ sends $\Lambda^2_{\pm}TM$ into
$\Lambda^2_{\mp}TM$, while ${\cal W}_{\pm}|\Lambda^2_{\mp}TM=0$.

A manifold $M$ is Einstein if and only if ${\cal B}=0$. It is
conformally flat if and only if ${\cal W}={\cal W}_{+}+{\cal W}_{-}$
vanishes. Recall also that $M$ is called self-dual (anti-self-dual)
if ${\cal W}_{-}=0$ (${\cal W}_{+}=0$).

For every $a\in\Lambda ^2TM$, define a skew-symmetric endomorphism
of $T_{\pi(a)}M$ by

\begin{equation}\label{cs}
g(K_{a}X,Y)=2g(a, X\wedge Y), \quad X,Y\in T_pM.
\end{equation}
If $\sigma\in\Lambda^2_{+}TM$ is a unit vector, then  $K_{\sigma}$
is a complex structure on the vector space $T_{\pi(\sigma)}$
compatible with the metric and the orientation of $M$. Conversely,
the $2$-vector $\sigma$ dual to one half of the K\"ahler $2$-form of
such a complex structure is a unit vector in $\Lambda^2_{+}TM$. Thus
the unit sphere subbunlde ${\cal Z}$ of $\Lambda^2_{+}TM$
parametrizes the complex structures on the tangent space of $M$
compatible with its metric and the orientation. This subbundle is
called the twistor space of $M$.

\smallskip

\noindent {\it Remark}. If we endow $\Lambda^2TM$ with the metric
$2g$, as many authors do, then the curvature operator acting on
$\Lambda^2TM$ is one half of the operator used here and the twistor
space of $M$ is the sphere subbundle of $\Lambda^2_{+}TM$ of radius
$\sqrt{2}$.

\smallskip

Let $(E_1,E_2,E_3,E_4)$ be a local oriented orthonormal frame of
$TM$. Set

\begin{equation}\label{s-basis}
\begin{array}{l}
s_1=E_1\wedge E_2+E_3\wedge E_4,  \hspace{1cm}   \bar{s}_1=E_1\wedge E_2-E_3\wedge E_4\\[6pt]
s_2=E_1\wedge E_3+E_4\wedge E_2, \hspace{1cm}    \bar{s}_2=E_1\wedge E_3-E_4\wedge E_2\\[6pt]
s_3=E_1\wedge E_4+E_2\wedge E_3,  \hspace{1cm}   \bar{s}_3=E_1\wedge
E_4-E_2\wedge E_3.
\end{array}
\end{equation}
Then $(s_1,s_2,s_3)$, resp. $(\bar{s}_1,\bar{s}_2,\bar{s}_3)$, is a local
oriented orthonormal frame of $\Lambda^2_{+}TM$, resp.
$\Lambda^2_{-}TM$.

For every $\sigma\in {\cal Z}$, the tangent space $T_{\pi(\sigma)}M$
has an orthonormal basis of the form $E',
K_{\sigma}E',E'',K_{\sigma}E''$. This basis yields the orientation
of $T_{\pi(\sigma)}$, so setting  $E_1=E', E_2=K_{\sigma}E',E_3=E'',
E_4=K_{\sigma}E''$ we obtain an oriented orthonormal basis for which
$\sigma=s_1$; for $E_1=E', E_2=E'', E_3=K_{\sigma}E',
E_4=-K_{\sigma}E''$, we have $\sigma=s_2$ and if $E_1=E', E_2=E'',
E_3=K_{\sigma}E'', E_4=K_{\sigma}E'$, we have $\sigma=s_3$.

\smallskip
The Levi-Civita connection $\nabla$ of $M$ induces a metric
connection on the bundle $\Lambda^2_{+}TM$ whose horizontal
distribution is tangent to the twistor space ${\cal Z}$. Thus we
have the decomposition $T{\cal Z}={\cal H}\oplus {\cal V}$ of the
tangent bundle of ${\cal Z}$ into horizontal and vertical
components. The vertical space at a point $\sigma\in{\cal Z}$ is the
space $${\cal V}_{\sigma}=\{V\in T_{\sigma}{\cal Z}:~
\pi_{\ast}V=0\}.$$ This is the tangent space to the fibre of ${\cal
Z}$ through $\sigma$, thus, considering $T_{\sigma}{\cal Z}$ as a
subspace of $T_{\sigma}(\Lambda^2_{+}TM)$ (as we shall always do),
${\cal V}_{\sigma}$ is the orthogonal complement of ${\Bbb R}\sigma$
in $\Lambda^2_{+}T_{\pi(\sigma)}M$. The map $V\ni{\cal
V}_{\sigma}\to K_{V}$ gives an identification of the vertical space
with the space of skew-symmetric endomorphisms of $T_{\pi(\sigma)}M$
that anti-commute with $K_{\sigma}$. Let $s$ be a local section of
${\cal Z}$ such that $s(p)=\sigma$ where $p=\pi(\sigma)$.
Considering $s$ as a section of $\Lambda^2_{+}TM$, we have
$\nabla_{X}s\in{\cal V}_{\sigma}$ for every $X\in T_pM$ since $s$
has a constant length. Moreover,
$X^h_{\sigma}=s_{\ast}X-\nabla_{X}s$ is the horizontal lift of $X$
at ${\sigma}$.

\smallskip

Denote by $\times$ the usual vector cross product on the oriented
$3$-dimensional vector space $\Lambda^2_{+}T_pM$, $p\in M$, endowed
with the metric $g$. Then it easy to check that
\begin{equation}\label{r-r}
g(R(a)b,c)=g({\cal R}(b\times c),a)
\end{equation}
for $a\in\Lambda^2T_pM$, $b,c\in\Lambda^2_{+}T_pM$, and
\begin{equation}\label{vp}
g(\sigma\times V,X\wedge K_{\sigma}Y)=g(\sigma\times V,K_{\sigma}X\wedge Y)=g(V,X\wedge Y)
\end{equation}
for $V\in{\cal V}_{\sigma}$, $X,Y\in T_pM$.

It is also easy to show that for every $a,b\in\Lambda^2_{+}T_pM$
\begin{equation}\label{com}
K_a\circ K_b=-g(a,b)Id + K_{a\times b}.
\end{equation}

The action of $SO(4)$ on $\Lambda^2{\Bbb R}^4$ preserves the
decomposition $\Lambda^2{\Bbb R}^4=\Lambda^2_{+}{\Bbb R}^4\oplus
\Lambda^2_{-}{\Bbb R}^4$. Thus, considering $S^2$ as the unit sphere
in $\Lambda^2_{+}{\Bbb R}^4$, we have an action of the group $SO(4)$
on $S^2$. Then, if $SO(M)$ denotes  the principal bundle of the
oriented orthonormal frames on $M$, the twistor space ${\cal Z}$ is
the associated bundle $SO(M)\times_{SO(4)} S^2$. It follows from the
Vilms theorem (see, for example, \cite[Theorem 9.59]{Besse}) that
the projection map $\pi:({\cal Z},h_t)\to (M,g)$ is a Riemannian
submersion with totally geodesic fibres (this can also be proved by
a direct computation).

Let $(G,x_1,...,x_4)$ be a local coordinate system of $M$ and let
$(E_1,...,E_4)$ be an oriented orthonormal frame of $TM$ on $G$. If
$(s_1,s_2,s_3)$ is the local frame of $\Lambda^2_{+}TM$ define by
(\ref{s-basis}), then $\widetilde x_{\alpha}=x_{\alpha}\circ\pi$,
$y_j(\sigma)=g(\sigma, (s_j\circ\pi)(\sigma))$, $1\leq \alpha\leq
4$, $1\leq j\leq 3$, are local coordinates of $\Lambda^2_{+}TM$ on
$\pi^{-1}(G)$.

   The horizontal lift $X^h$ on $\pi^{-1}(G)$ of a vector field
$$
X=\sum_{\alpha=1}^4 X^{\alpha}\frac{\partial}{\partial x_{\alpha}}
$$
is given by
\begin{equation}\label{hl}
X^h=\sum_{\alpha=1}^4 (X^{\alpha}\circ\pi)\frac{\partial}{\partial
\widetilde{x}_{\alpha}}
-\sum_{j,k=1}^3y_j(g(\nabla_{X}s_j,s_k)\circ\pi)\frac{\partial}{\partial
y_k}.
\end{equation}
Hence
\begin{equation}\label{Lie-1}
[X^h,Y^h]=[X,Y]^h+\sum_{j,k=1}^3y_j(g(R(X\wedge
Y)s_j,s_k)\circ\pi)\frac{\partial}{\partial y_k}.
\end{equation}
for every vector fields $X,Y$ on $G$. Let $\sigma\in{\cal Z}$. Using
the standard identification $T_{\sigma}(\Lambda^2_{+}T_pM)\cong
\Lambda^2_{+}T_pM$ formula (\ref{Lie-1}) can be rewritten as
\begin{equation}\label{Lie-2}
[X^h,Y^h]_{\sigma}=[X,Y]^h_{\sigma}+R_{p}(X\wedge Y)\sigma, \quad
p=\pi(\sigma).
\end{equation}

Denote by $D$ the Levi-Civita connection of $({\cal Z},h_t)$. Then we have the following.
\begin{lemma}\label{LC} \rm ({\cite{DM})}
If $X,Y$ are vector fields on $M$ and $V$ is a vertical vector field on ${\cal Z}$, then
\begin{equation}\label{D-hh}
(D_{X^h}Y^h)_{\sigma}=(\nabla_{X}Y)^h_{\sigma}+\frac{1}{2}R_{p}(X\wedge
Y)\sigma,
\end{equation}
\begin{equation}\label{D-vh}
(D_{V}X^h)_{\sigma}={\cal H}(D_{X^h}V)_{\sigma}=-\frac{t}{2}(R_{p}(\sigma\times V)X)^h_{\sigma}
\end{equation}
where $\sigma\in{\cal Z}$, $p=\pi(\sigma)$ and ${\cal H}$ means "the
horizontal component".
\end{lemma}
{\bf Proof}. Identity (\ref{D-hh}) follows from the Koszul formula
for the Levi-Civita connection and (\ref{Lie-2}).

Let $W$ be e vertical vector field on ${\cal Z}$. Then
$$
h_t(D_{V}X^h,W)=-h_t(X^h,D_{V}W)=0
$$
since the fibres are totally geodesic submanifolds, so $D_{V}W$ is a
vertical vector field. Therefor $D_{V}X^h$ is a horizontal vector
field. Moreover, $[V,X^h]$ is a vertical vector field, hence
$D_{V}X^h={\cal H}D_{X^h}V$. Then
$$
h_t(D_{V}X^h,Y^h)=h_t(D_{X^h}V,Y^h)=-h_t(V,D_{X^h}Y^h).
$$
Now (\ref{D-vh}) follows from (\ref{D-hh}) and (\ref{r-r}). {\it
Q.E.D.}

\section{Compatible almost complex structures on twistor spaces}

Let $f:{\cal Z}\to {\cal Z}$ be  a morphism of the bundle ${\cal Z}$, i.e. a smooth map with
$\pi\circ f=\pi$. Following \cite{De} we define an almost complex structure ${\cal J}_f$ on the
$6$-manifold ${\cal Z}$ setting
$$
{\cal J}_fV=\sigma\times V ~ {\mbox~for} ~ V\in{\cal V}_{\sigma}
$$
$$
{\cal J}_fX^h_{\sigma}=(K_{f(\sigma)}X)^h_{\sigma}~ {\mbox~for} ~
X\in T_{\pi(\sigma)}M.
$$
Note that, since the fibres of ${\cal Z}$ are spheres, the
restriction of ${\cal J}_f$ to any fibre is the standard complex
structure of the unite sphere.

In the case when $f=Id$, the almost complex structure ${\cal J}_f$ coincides with that defined by
Atiyah-Hitchin-Singer \cite{AHS}. In this case the almost complex structure ${\cal J}_f$ is
integrable if and only if the base manifold $M$ is anti-self-dual \cite{AHS}. If $f$ is the
antipodal map $f(\sigma)=-\sigma$, ${\cal J}_f$ is the conjugate structure of the almost complex
structure defined by Eells-Salamon \cite{ES}. This structure is never integrable \cite{ES}.

\smallskip
\noindent {\it Remark}. One can also consider the sphere bundle ${\cal Z}^{-}$ in $\Lambda^2_{-}TM$
as the twistor space of $M$. Then the Atiyah-Hitchin-Singer almost complex structure is defined as
${\cal J}V=-\sigma\times V$ for $V\in{\cal V}_{\sigma}$ and ${\cal
J}X^h_{\sigma}=(K_{\sigma}X)^h_{\sigma}$ for $X\in T_{\pi(\sigma)}M$. It is integrable if and only
if $M$ is self-dual. The complex projective space ${\Bbb C}{\Bbb P}^2$ with the Fubini-Study metric
is self-dual but not anti-self-dual. Thus the Atiyah-Hitchin-Singer almost complex structure on
${\cal Z}$ is not integrable while it is integrable on ${\cal Z}^{-}$. This is one of the reasons
some authors to consider the sphere bundle ${\cal Z}^{-}$ in $\Lambda^2_{-}TM$ instead of that in
$\Lambda^2_{+}TM$.

\smallskip

The almost complex structure ${\cal J}_f$ is compatible with the Riemannian metrics $h_t, t>0$
defined above and let $\Omega(A,B)=h_t({\cal J}_fA,B)$ be the K\"ahler $2$-form of the almost
Hermitian manifold $(Z, h_t, {\cal J}_f)$. We now compute the covariant derivative of $\Omega$.

\begin{prop}\label{D-Om}
Let $\sigma\in{\cal Z}$, $X,Y,Z\in T_pM$, $p=\pi(\sigma)$,
$U,V,W\in{\cal V}_{\sigma}$. Then
\begin{equation}\label{DOm-3hor}
(D_{X^h_{\sigma}}\Omega)(Y^h_{\sigma},Z^h_{\sigma})=2g({\cal
V}f_{\ast}(X^h_{\sigma}),Y_\wedge Z);
\end{equation}
\begin{equation}\label{DOm-2hor-ver}
(D_{X^h_{\sigma}}\Omega)(Y^h_{\sigma},U)=-\frac{t}{2}g({\cal R}(U),X\wedge Y) +\frac{t}{2}g({\cal
R}(\sigma\times U), X\wedge K_{f(\sigma)}Y);
\end{equation}
\begin{equation}\label{DOm-ver-2hor}
(D_{U}\Omega)(Y^h_{\sigma},Z^h_{\sigma})=-\frac{t}{2}g({\cal
R}(\sigma\times U),Y\wedge K_{f(\sigma)}Z+K_{f(\sigma)}Y\wedge Z)
+2g(f_{\ast}U,Y\wedge Z);
\end{equation}
\begin{equation}\label{DOm-hor-ver}
(D_{X^h_{\sigma}}\Omega)(U,V)=0,\quad
(D_{U}\Omega)(Y^h_{\sigma},V)=0,\quad (D_{U}\Omega)(V,W)=0.
\end{equation}
\end{prop}
{\bf Proof}.  Extend the tangent vector $Y, Z$ to vector fields in a
neighbourhood of $p$ such that $\nabla Y|_p=\nabla Z|_p=0$.

To prove the first formula, we note that
$$
(D_{X^h_{\sigma}}\Omega)(Y^h_{\sigma},Z^h_{\sigma})=X^h_{\sigma}(h_t({\cal J}_fY^h,Z^h))-h_t({\cal
J}_fD_{X^h_{\sigma}}Y^h,Z^h)+ h_t(Y^h,{\cal J}_fD_{X^h_{\sigma}}Z^h).
$$
The vectors ${\cal J}_fD_{X^h_{\sigma}}Y^h$ and ${\cal J}_fD_{X^h_{\sigma}}Z^h$  are vertical in view of (\ref{D-hh}). Hence \\
$(D_{X^h_{\sigma}}\Omega)(Y^h_{\sigma},Z^h_{\sigma})=X^h_{\sigma}(h_t({\cal J}_fY^h,Z^h))$.  Let
$s$ be a section of the bundle ${\cal Z}$ around $p$ such that $s(p)=\sigma$ and $\nabla s|_p=0$.
Then
$$
\begin{array}{c}
X^h_{\sigma}(h_t({\cal J}_fY^h,Z^h))=X(h_t({\cal J}_fY^h,Z^h))\circ s)=X(g(K_{f(s)}Y,Z))=\\[6pt]
2X(g(f\circ s,Y\wedge Z))=2g(\nabla_{X}(f\circ s),Y\wedge Z).
\end{array}
$$
The map $\widetilde s=f\circ s$ is a section of ${\cal Z}$ with
$\widetilde s(p)=f(\sigma)$. Then
$$
X^h_{f(\sigma)}+\nabla_{X}(f\circ s)=\widetilde
s_{\ast}(X)=f_{\ast\,\sigma}(s_{\ast\,p}(X))=f_{\ast}(X^h_{\sigma}).
$$
Therefore $\nabla_{X}(f\circ s)={\cal V}f_{\ast}(X^h_{\sigma})$.
Thus  $(D_{X^h_{\sigma}}\Omega)(Y^h_{\sigma},Z^h_{\sigma})= 2g({\cal
V}f_{\ast}(X^h_{\sigma}),Y_\wedge Z)$.

Extend $U$ to a vertical vector field in a neighbourhood of
$\sigma$.  Identities  (\ref{D-hh}) and (\ref{D-vh}) imply that
$$
\begin{array}{c}
(D_{X^h_{\sigma}}\Omega)(Y^h_{\sigma},U)=-h_t({\cal J}_fD_{X^h_{\sigma}}Y^h,U)-h_t({\cal J}_fY^h_{\sigma},D_{X^h_{\sigma}}U)=\\[6pt]
-\displaystyle{\frac{t}{2}}g(\sigma\times
R(X,Y)\sigma,U)+\displaystyle{\frac{t}{2}}g(K_{f(\sigma)}Y,R(\sigma\times U)X)
\end{array}
$$
This gives the second formula of the lemma since, in view of
(\ref{r-r}),
$$
g(\sigma\times R(X,Y)\sigma,U)=-g(R(X,Y)\sigma,\sigma\times U)=g({\cal R}(U),X\wedge Y).
$$

Next, we have
\begin{equation}\label{v-hh}
(D_{U}\Omega)(Y^h_{\sigma},Z^h_{\sigma})=U(h_t({\cal J}_fY^h,Z^h))+h_t(D_{U}Y^h,{\cal
J}_fZ^h)-h_t({\cal J}_fY^h,D_{U}Z^h).
\end{equation}
Moreover, $f=\sum_{i=1}^3 (y_i\circ f)(s_i\circ\pi)$, therefore
$$
U(h_t({\cal J}_fY^h,Z^h))=2\sum_{i=1}^3U((y_i\circ f)g(s_i,Y\wedge
Z)\circ\pi)=2\sum_{i=1}^3U(y_i\circ f)g(s_i,Y\wedge Z)_p.
$$
The map $f$ sends fibres to fibres, hence $f_{\ast}$ sends vertical vectors to vertical vectors. In
particular, $f_{\ast}U=\sum_{i=1}^3 U(y_i\circ f)\displaystyle{(\frac{\partial}{\partial
y_i}})_{f(\sigma)}$. It follows that $U(h_t({\cal J}_fY^h,Z^h))=2g(f_{\ast}U,Y\wedge Z)$ and the
third formula of the lemma follows from (\ref{v-hh}) and (\ref{D-vh}).

To prove the remaining formulas fix a point $\sigma\in{\cal Z}$ and
set $p=\pi(\sigma)$. Take an oriented orthonormal frame
$(E_1,...,E_4)$ of $M$ around the point $p$ such that $\nabla
E_{\alpha}|_p=0$, $\alpha=1,...,4$, and define an oriented
orthonormal frame $(s_1,s_2,s_3)$ of $\Lambda^2_{+}TM$ by means of
(\ref{s-basis}). We have $\nabla s_i|_p=0$, $i=1,2,3$,  for the
latter frame. Choose also a local coordinate system $(x_1,...,x_4)$
of $M$ near $p$, then define local coordinates $(x_{\alpha},y_i)$
$\alpha=1,...,4$, $i=1,2,3$, on ${\cal Z}$ as above.

Every section $a$ of $\Lambda^2_{+}TM$ on an open set $G$ gives rise
to a vertical vector field $\widetilde a$ on $\pi^{-1}(G)$ defined
by
$$
{\widetilde
a}_{\tau}=a\circ\pi(\tau)-g(a\circ\pi(\tau),\tau)\tau,\quad
\tau\in\pi^{-1}(G).
$$
Note that, around every point of ${\cal Z}$, there exists a frame of
vertical vector fields of this type.

Further on, we shall use this notation without explicitly saying so.

Now take sections $a$ and $b$  of ${\cal Z}$ defined in a
neighbourhood of $p=\pi(\sigma)$ and such that $a(p)=U$, $b(p)=V$,
$\nabla a|_p=\nabla b|_p=0$. Let $\widetilde a$ and $\widetilde b$
be the vertical vector fields associated to $a$ and $b$. Then
$\widetilde a_{\sigma}=U$, $\widetilde b_{\sigma}=V$ and
\begin{equation}\label{h-vv}
(D_{X^h_{\sigma}}\Omega)(U,V)=X^h_{\sigma}(h_t({\cal J}_f\widetilde a,\widetilde b))-h_t({\cal
J}_fD_{X^h_{\sigma}}\widetilde a,V)+ h_t(U,{\cal J}_fD_{X^h_{\sigma}}\widetilde b).
\end{equation}
Set
$$
\widetilde a=\sum_{i=1}^3 \widetilde{a}_i\frac{\partial}{\partial
y_i},\quad \widetilde b=\sum_{i=1}^3
\widetilde{b}_i\frac{\partial}{\partial y_i}.
$$
Then
\begin{equation}\label{ai}
\widetilde a_i=\sum_{j=1}^3(\delta_{ij}-y_iy_j)(g(a,s_j)\circ\pi),
\end{equation}
and similar for $\widetilde b_i$. Moreover,
\begin{equation}\label{ja}
{\cal J}_f\widetilde a=(y_2\widetilde{a}_3-y_3\widetilde{a}_2)\frac{\partial}{\partial y_1}
+(y_3\widetilde{a}_1-y_1\widetilde{a}_3)\frac{\partial}{\partial y_2}
+(y_1\widetilde{a}_2-y_2\widetilde{a}_1)\frac{\partial}{\partial y_3}.
\end{equation}
Hence
$$
h_t({\cal J}_f\widetilde a,\widetilde b))=(y_2\widetilde{a}_3-y_3\widetilde{a}_2)\widetilde{b}_1
+(y_3\widetilde{a}_1-y_1\widetilde{a}_3)\widetilde{b}_2
+(y_1\widetilde{a}_2-y_2\widetilde{a}_1)\widetilde{b}_3.
$$
If $X=\sum_{\alpha=1}^4
X^{\alpha}(\displaystyle{\frac{\partial}{\partial x_{\alpha}}})_p$,
we have $X^h_{\sigma}=\sum_{\alpha=1}^4
X^{\alpha}(\displaystyle{\frac{\partial}{\partial
\widetilde{x}_{\alpha}}})_{\sigma}$, hence
$$
X^h_{\sigma}(\widetilde{a}_i)=\sum_{j=1}^3(\delta_{ij}-y_iy_j)X(g(a,s_j))=0
$$
since $\nabla_{X}a=\nabla_{X}s_j=0$. Similarly, $X^h_{\sigma}(\widetilde{b}_i)=0$, $i=1,2,3$. It
follows that $X^h_{\sigma}(h_t({\cal J}_f\widetilde a,\widetilde b))=0$. Using (\ref{hl}) and
(\ref{ai}), one obtains  by a straightforward computation that $[X^h,\widetilde
a]_{\sigma}=\widetilde{(\nabla_X a)}_{\sigma}=0$. Hence $D_{X^h_{\sigma}}\widetilde a=
-D_{{\widetilde a}_{\sigma}}X^h\in{\cal H}_{\sigma}$ in view of (\ref{D-vh}). Then $h_t({\cal
J}_fD_{X^h_{\sigma}}\widetilde a,V)=0$. Similarly $h_t(U,{\cal J}_fD_{X^h_{\sigma}}\widetilde
b)=0$. Thus, $(D_{X^h_{\sigma}}\Omega)(U,V)=0$ by  (\ref{h-vv}). Also
$$
(D_{U}\Omega)(Y^h_{\sigma},V)=U(h_t({\cal J}_fY^h,\widetilde b))-h_t({\cal
J}_fD_{U}Y^h,V)-h_t({\cal J}_fY^h,D_{U}\widetilde b)=0
$$
since ${\cal J}_fY^h$, ${\cal J}_fD_{U}Y^h$ are horizontal vectors and $D_{U}\widetilde b$ is
vertical.

Finally, the identity $(D_{U}\Omega)(V,W)=0$ is a consequence of the fact that the fibres of ${\cal
Z}$ are totally geodesic submanifolds and ${\cal J}_f$ preserves the vertical distribution.

Proposition~\ref{D-Om} and the formula $$d\Omega(A,B,C)=
\sum_{cyc}(D_{A}\Omega)(B,C)$$ give the following

\begin{cor}\label{diff-Om}
Let $\sigma\in{\cal Z}$, $X,Y,Z\in T_pM$, $p=\pi(\sigma)$,
$U,V,W\in{\cal V}_{\sigma}$. Then
$$
d\Omega(X^h,Y^h,Z^h)=2g({\cal V}f_{\ast}(X^h),Y\wedge Z)+2g({\cal
V}f_{\ast}(Y^h),Z\wedge X)+2g({\cal V}f_{\ast}(Z^h),X\wedge Y)
$$
$$
d\Omega(X^h,Y^h,U)=g(2f_{\ast}U-t{\cal R}(U),X\wedge Y)
$$
$$
d\Omega(X^h,U,V)=0,\quad d\Omega(U,V,W)=0.
$$
\end{cor}

\begin{cor}\label{delta-Om}
Let $\sigma\in{\cal Z}$, $X\in T_pM$, $p=\pi(\sigma)$, $U\in{\cal
V}_{\sigma}$. Then
\begin{equation}\label{delta-hor}
\begin{array}{c}
(\delta\Omega)(X^h_{\sigma})=Trace\{T_{p}M\ni A\to 2g({\cal V}f_{\ast}(A^h_{\sigma}),X\wedge A)\}=\\[6pt]
Trace\{{\cal V}_{\sigma}\ni\tau\to g({\cal
V}f_{\ast}((K_{\tau}X)^h_{\sigma}),\tau)\}.
\end{array}
\end{equation}
\begin{equation}\label{delta-ver}
\delta\Omega(U)=-tg({\cal R}(\sigma\times U),f(\sigma)).
\end{equation}
\end{cor}
{\bf Proof}. Let $E_1,...,E_4$ be an orthonormal basis of $T_pM$,
$p=\pi(\sigma)$ and $\tau_1,\tau_2$ a $g$-orthonormal basis of
${\cal V}_{\sigma}$. Then, by (\ref{DOm-3hor}) and
(\ref{DOm-hor-ver}),
$$
\begin{array}{l}
\delta\Omega(X^h_{\sigma})=-\sum_{i=1}^4(D_{(E_i)^h_{\sigma}}\Omega)((E_i)^h_{\sigma},X^h_{\sigma})-\sum_{m=1}^2(D_{\tau_m}\Omega)(\tau_m,X^h_{\sigma})=\\[6pt]
-2\sum_{i=1}^4g({\cal V}f_{\ast}((E_i)^h_{\sigma}),E_i\wedge X)=-2\sum_{i=1}^4\sum_{m=1}^2g(f_{\ast}((E_i)^h_{\sigma}),\tau_m)g(\tau_m,E_i\wedge X)=\\[6pt]
\sum_{i=1}^2\sum_{m=1}^2g(f_{\ast}((E_i)^h_{\sigma}),\tau_m)g(K_{\tau_m}X,E_i)=\sum_{m=1}^2g(f_{\ast}((K_{\tau_m}X)^h_{\sigma}),\tau_m).
\end{array}
$$
This proves (\ref{delta-hor}).

In view of ({\ref{DOm-2hor-ver}) and (\ref{DOm-hor-ver}), we have
$$
\delta\Omega(U)=-\frac{t}{2}g({\cal R}(\sigma\times U),\sum_{i=1}^4E_i\wedge K_{f(\sigma)}E_i)
$$
Moreover, for $Y,Z\in T_pM$,
$$
\begin{array}{c}
\sum_{i=1}^4g(E_i\wedge K_{f(\sigma)}E_i,Y\wedge Z)=\\[6pt]
\frac{1}{2}\sum_{i=1}^4[-g(Y,E_i)g(K_{f(\sigma)}Z,E_i)+g(Z,E_i)g(K_{f(\sigma)}Y,E_i)]=\\[6pt]
g(K_{f(\sigma)}Y,Z)=2g(f(\sigma),Y\wedge Z).
\end{array}
$$
Thus $\sum_{i=1}^4E_i\wedge K_{f(\sigma)}E_i=2f(\sigma)$ and the
second formula of the corollary is proved.
\smallskip

Denote the Nijenhuis tensor of ${\cal J}_f$ by $N$. The next statement follows from
Proposition~\ref{D-Om}, identity (\ref{vp}) and the well-known formula
$$
h_t(N(A,B),C)=(D_A\Omega)({\cal J}_fB,C)-(D_{{\cal J}_fB}\Omega)(A,C)-(D_B\Omega)({\cal
J}_fA,C)+(D_{{\cal J}_fA}\Omega)(B,C).
$$
\begin{cor}\label{N}
Let $\sigma\in{\cal Z}$, $X,Y,Z\in T_pM$, $p=\pi(\sigma)$,
$U,V\in{\cal V}_{\sigma}$. Then
$$
\begin{array}{c}
h_t(N(X^h_{\sigma},Y^h_{\sigma}),Z^h_{\sigma})=2g({\cal
V}f_{\ast}(X^h_{\sigma}),K_{f(\sigma)}Y\wedge Z)-
2g({\cal V}f_{\ast}(Y^h_{\sigma}),K_{f(\sigma)}X\wedge Z)\\[6pt]
+2g({\cal V}f_{\ast}((K_{f(\sigma)}X)^h_{\sigma}),Y\wedge Z)-2g({\cal V}f_{\ast}((K_{f(\sigma)}Y)^h_{\sigma}),X\wedge Z)\\[10pt]

h_t(N(X^h_{\sigma},Y^h_{\sigma}),U)=-tg({\cal R}(X\wedge K_{f(\sigma)}Y+K_{f(\sigma)}X\wedge Y),U)\\[6pt]
-tg({\cal R}(X\wedge Y-K_{f(\sigma)}X\wedge K_{f(\sigma)}Y) ,\sigma\times U)\\[10pt]

h_t(N(X^h_{\sigma},U),Z^h_{\sigma})=-2g(f(\sigma)\times f_{\ast}(U),X\wedge Z)+2g(f_{\ast}(\sigma\times U),X\wedge Z)\\[10pt]

h_t(N(X^h_{\sigma},U),V)=0\quad N(U,V)=0.
\end{array}
$$
\end{cor}

Since $h_t(N(X^h,U),Z^h)=-2g({\cal J}_{f}f_{\ast}(U),X\wedge Z)+2g(f_{\ast}({\cal J}_fU),X\wedge
Z)$, we have the following.

\begin{cor}\label{hol f} \rm ({\cite{De})}
${\cal H}(N(X^h,U))=0$ if and only if the restriction of $f$ to
every fibre is a holomorphic map.
\end{cor}

\smallskip

\section{Gray-Hervella classes of the almost complex structures  ${\cal J}_{\omega}$ }

In what follows we use the same notation for the Gray-Hervella
classes as in \cite{GH}. For example, ${\cal K}$ is the class of
K\"ahler manifolds, ${\cal W}_1$ is the class of nearly K\"ahler
manifolds, ${\cal W}_2$ is the class of almost K\"ahler manifolds,
${\cal W}_3\oplus {\cal W}_4$ is the class of Hermitian manifolds,
${\cal W}_1\oplus {\cal W}_2\oplus {\cal W}_3$ is the class of
semi- K\"ahler or balanced manifolds, etc.

Let $(g,J)$ be an almost Hermitian structure on a four-manifold $M$.
Define a section $\omega$ of $\Lambda^2TM$ by
$$
g(\omega,X\wedge Y)=\frac{1}{2}g(JX,Y),~X,Y\in TM.
$$
Clearly, at any point, $\omega$ is the dual $2$-vector of one half
of the K\"ahler $2$-form $F$ of the almost Hermitian manifold
$(M,g,J)$. Consider $M$ with the orientation yielded by the almost
complex structure $J$. Then $\omega$ is a section of the twistor
bundle ${\cal Z}$. As in \cite{De}, define a bundle map $f:{\cal
Z}\to {\cal Z}$ setting $ f=\omega\circ\pi.$ Since the restriction
of $f$ to any fibre is a constant map, $f_{\ast}|{\cal V}=0.$ We
also have
\begin{equation}\label{hor}
f_{\ast}(X^h_{\sigma})=X^h_{\omega(p)}+\nabla_X\omega,
\end{equation}

\noindent where $p=\pi(\sigma)$ and $X\in T_pM$. Note that
$$2g(\nabla_X\omega,Y\wedge Z)=(\nabla_XF)(Y,Z).$$

\smallskip

Denote by ${\cal J}_{\omega}$ the almost complex structure on ${\cal
Z}$ determined by the map $f$ defined by $\omega$. In the next
theorem we determine the Gray-Hervella classes of the almost
Hermitian manifolds  $({\cal Z},h_t, {\cal J}_{\omega})$.

\smallskip

\begin{th}\label{GH}
Let $(M,g,J)$ be an almost Hermitian $4$-manifold with K\"ahler
$2$-vector $\omega$, self-dual Weyl tensor $W_+$ and scalar
curvature $s$. The possible Gray-Hervella classes of its twistor
space $({\cal Z},h_t, {\cal J}_{\omega})$ are ${\cal W},~ {\cal K},
~{\cal W}_3, ~{\cal H}={\cal W}_3\oplus {\cal W}_4, ~{\cal SK}={\cal
W}_1\oplus {\cal W}_2\oplus {\cal W}_3, ~{\cal G}_1={\cal W}_1\oplus
{\cal W}_3\oplus {\cal W}_4$ and ${\cal G}_2={\cal W}_2\oplus {\cal
W}_3\oplus {\cal W}_4$. Moreover

\smallskip

\noindent $(i)$ $({\cal Z},h_t, {\cal J}_{\omega})\in {\cal K}$ if
and only if $(M,g,J)$ is  K\"ahler and Ricci flat.

\medskip

\noindent $(ii)$ $({\cal Z},h_t, {\cal J}_{\omega})\in {\cal SK}\cap
{\cal H}={\cal W}_3$ if and only if $(M,g,J)$ is K\"ahler and scalar
flat.

\medskip

\noindent $(iii)$ $({\cal Z},h_t, {\cal J}_{\omega})\in {\cal
H}={\cal W}_3\oplus {\cal W}_4$ if and only if $(M,g,J)$ is
Hermitian and
$${\cal
W}_{+}(\sigma)=\displaystyle{\frac{s}{2}g(\sigma,\omega)\omega-\frac{s}{6}\sigma}$$
for all $\sigma\in\Lambda^2_{+}TM$.

\medskip

\noindent $(iv)$ $({\cal Z},h_t, {\cal J}_{\omega})\in{\cal
SK}={\cal W}_1\oplus {\cal W}_2\oplus {\cal W}_3$ if and only if
$(M,g,J)$ is almost K\"ahler and $${\cal
W}_{+}(\omega)=-\displaystyle{\frac{s}{6}}\omega.$$

\medskip

\noindent $(v)$ $({\cal Z},h_t, {\cal J}_{\omega})\in{\cal
G}_1={\cal W}_1\oplus {\cal W}_3\oplus {\cal W}_4$ if and only if
$(M,g,J)$ is Hermitian.

\medskip

\noindent $(vi)$ $({\cal Z},h_t, {\cal J}_{\omega})\in{\cal
G}_2={\cal W}_2\oplus {\cal W}_3\oplus {\cal W}_4$ if and only if
$${\cal
W}_{+}(\sigma)=\displaystyle{\frac{s}{2}g(\sigma,\omega)\omega-\frac{s}{6}\sigma}$$
for all $\sigma\in\Lambda^2_{+}TM$.

\end{th}

{\bf Proof}. To determine the possible Gray-Hervella classes of the twistor space $({\cal Z},h_t,
{\cal J}_{\omega})$  we shall need several technical lemmas.

Given a point ${\cal Z}$, we take a basis $E_1,E_2=JE_1,E_3,E_4=JE_3$ of $T_{\pi(\sigma)}M$. Such a
basis induces the orientation of $M$ we have chosen and we define $s_1,s_2,s_3$ and
$\bar{s}_1,\bar{s}_2,\bar{s}_3$ via (\ref{s-basis}) (so $\omega=s_1$). This notation will be used in the proofs
of the next statements.

\begin{lemma}\label{K} $({\cal Z},h_t, {\cal J}_{\omega})\in {\cal K}$ if and only if
$(M,g,J)$ is K\"ahler and Ricci flat.
\end{lemma}
{\bf Proof}. It follows from  Proposition~\ref{D-Om} and (\ref {hor})  that $({\cal Z},h_t, {\cal
J}_{\omega})$ is K\"ahler if and only if $(M,g,J)$ is K\"ahler and for every $\sigma\in{\cal Z}$,
$U\in{\cal V}_{\sigma}$ and $X,Y\in T_{\pi(\sigma)}M$

\smallskip

\noindent $(i)$ $-g({\cal R}(U),X\wedge Y)+g({\cal R}(\sigma\times
U),X\wedge JY)=0,$

\medskip

\noindent $(ii)$ $g({\cal R}(U),X\wedge JY + JX\wedge Y)=0.$

\smallskip

The latter identity implies
\begin{equation}\label{23}
g({\cal R}(U),s_2)=g({\cal R}(U),s_3)=0.
\end{equation}
It follows from identity $(i)$ that
\begin{equation}\label{1234}
g({\cal R}(U),E_1\wedge E_2)=g({\cal R}(U),E_3\wedge E_4)=0,
\end{equation}
\begin{equation}\label{1235}
\begin{array}{c}
g({\cal R}(U),E_1\wedge E_3)=g({\cal R}(\sigma\times U),E_1\wedge
E_4),
\\[6pt] g({\cal R}(U),E_3\wedge E_1)=g({\cal R}(\sigma\times U),E_3\wedge E_2).
\end{array}
\end{equation}
We obtain from (\ref{1234}) that $g({\cal R}(U),\bar{s}_1)=g({\cal R}(U),s_1)=0$. Thus $g({\cal
R}(U),s_i)=0$ for $i=1,2,3$ and every $U\in \Lambda^2_{+}T_pM$. It follows from (\ref{1235}) that
$$
g({\cal R}(\sigma\times U),\bar{s}_3)=0.
$$
Moreover, identities (\ref{23}) and (\ref{1235}) imply
$$
g({\cal R}(U),\bar{s}_2)=2g({\cal R}(U),E_1\wedge E_3)=2g({\cal R}(\sigma\times U),E_1\wedge
E_4)=g({\cal R}(\sigma\times U),\bar{s}_3).
$$
Therefore $g({\cal R}(U),\bar{s}_2)=g({\cal R}(U),\bar{s}_3)=0$, thus
$g({\cal R}(U),\bar{s}_i)=0$, $i=1,2,3$. Hence ${\cal R}(U)=0$ for every
$U\in \Lambda^2_{+}T_pM$. This shows that if $({\cal Z},h_t, {\cal
J}_{\omega})$ is K\"ahler, then $(M,g,J)$ is a K\"ahler and Ricci
flat manifold.

Conversely, suppose that $(M,g,J)$ is K\"ahler and Ricci flat.  Using the curvature decomposition
(\ref{cur-dec}), the K\"ahler curvature identities and the first Bianchi identity, one can  see
that
$$
g({\cal R}(s_1),s_1)=\frac{s}{3}s_1,\quad {\cal R}(s_2)={\cal R}(s_3)=0.
$$
This implies the well-known fact (which can be traced back to \cite{Ga}) that the eigenvalues of
the operator ${\cal W}_{+}$ on a K\"ahler surface are
$\displaystyle{\frac{s}{3}},-\frac{s}{6},-\frac{s}{6},$. It follows that ${\cal R}(U)=0$ for every
$U\in \Lambda^2_{+}T_pM$, thus identities $(i)$ and $(ii)$ obviously are satisfied.

\begin{lemma}\label{W124} $({\cal Z},h_t, {\cal J}_{\omega})\in {\cal W}_1\oplus {\cal
W}_2\oplus {\cal W}_4$ if and only if  $({\cal Z},h_t, {\cal
J}_{\omega})\in {\cal K}$.
\end{lemma}
{\bf Proof}. The condition for $({\cal Z},h_t, {\cal J}_{\omega})$
to be in the class ${\cal W}_1\oplus {\cal W}_2\oplus {\cal W}_4$
is
\begin{equation}\label{124}
\begin{array}{c}
(D_A\Omega)(B,C)+(D_{{\cal J}_{\omega}A}\Omega)(\cal J_{\omega}B,C)=\\[6pt]
-\frac{1}{2}\{h_t(A,B)\delta\Omega(C)
-h_t(A,C)\delta\Omega(B)-h_t(A,{\cal J}_{\omega}B)\delta\Omega({\cal J}_{\omega}C)\\[6pt]
+h_t(A,{\cal J}_{\omega}C)\delta\Omega({\cal J}_{\omega}B)\}
\end{array}
\end{equation}
for every $A,B,C\in T{\cal Z}$. Proposition~\ref{D-Om} and (\ref
{hor}) imply that this condition is satisfied if and only if for
every $\sigma\in{\cal Z}$, $X,Y,Z\in T_{\pi(\sigma)}M$ and
$U,V,W\in{\cal V}_{\sigma}$ we have
\smallskip

\noindent $(i)$ $\begin{array}{c}
(\nabla_X F)(Y,Z)+(\nabla_{JX}
F)(JY,Z)=\\-\frac{1}{2}\{g(X,Y)\delta F(Z)-g(X,Z)\delta
F(Y)-g(X,JY)\delta F(JZ)+g(X,JZ)\delta F(JY)\},
\end{array}$

\medskip

\noindent $(ii)$ $\begin{array}{c}
-g({\cal R}(U),X\wedge Y+JX\wedge JY)+g({\cal R}(\sigma\times U),X\wedge JY-JX\wedge Y)=\\
g(X,Y)g({\cal R}(\sigma\times U),\omega)+g(X,JY)g({\cal R}(U),\omega),
\end{array}$

\medskip

\noindent $(iii)$ $g({\cal R}(U),X\wedge Y-JX\wedge JY)+g({\cal R}(\sigma\times U),X\wedge
JY+JX\wedge Y)=0$,

\medskip

\noindent $(iv)$ $g(U,V)\delta F(X)-g(U,\sigma\times V)\delta F(JX)=0$,

\medskip

\noindent $(v)$ $\begin{array}{c} g(U,V)g({\cal R}(\sigma\times W),\omega)-g(U,W)g({\cal R}(\sigma\times V),\omega)+\\[6pt]
g(U,\sigma\times V)g({\cal R}(W),\omega)-g(U,\sigma\times W)g({\cal
R}(V),\omega)=0.\end{array}$

\medskip

Clearly, identity $(v)$, obtained from (\ref{124}) for vertical vectors $A=U$, $B=V$, $C=W$, holds
when $U=0$. If $U\neq 0$, then $U,\sigma\times U$ is a basis of ${\cal V}_{\sigma}$ and it is easy
to check that this identity is also satisfied. Thus identity $(v)$  does not impose any restriction
on the base manifold $M$. Identity $(iv)$ implies that $\delta F=0$. Then it follows from $(i)$
that $(\nabla_X F)(Y,Z)+(\nabla_{JX} F)(JY,Z)=0$. It is well-known (and easy to see) that, in
dimension $4$, the latter identity is equivalent to $dF=0$. Take a point $p\in M$ and let $X\in
T_pM$ be a unit vector. For every point $\sigma\in{\cal Z}$ with $\pi(\sigma)=p$ and every $U\in
{\cal V}_{\sigma}$, identity $(ii)$ with $Y=JX$ gives $2g({\cal R}(U),X\wedge JX)=g({\cal
R}(U),s_1)$. Thus we have  $2g({\cal R}(U),E_1\wedge E_2)=g({\cal R}(U),s_1)$ and $2g({\cal
R}(U),E_3\wedge E_4)=g({\cal R}(U),s_1)$. This implies
$$
g({\cal R}(U),E_1\wedge E_2)=g({\cal R}(U),E_3\wedge E_4)=g({\cal R}(U),s_1)=0
$$
since $s_1=E_1\wedge E_2+E_3\wedge E_4$. It follows that
\begin{equation}\label{i1}
g({\cal R}(U),s^{-}_1)=0,\quad g({\cal R}(s_1),s_1)=g({\cal R}(s_2),s_1)=g({\cal R}(s_3),s_1)=0.
\end{equation}
Identity $(iii)$ with $X=E_1$, $Y=E_3$ becomes
$$
g({\cal R}(U),s_2)+g({\cal R}(\sigma\times U),s_3)=0,\quad U\in{\cal
V}_{\sigma}.
$$
Applying the latter identity for $\sigma=s_2$ and $\sigma=s_3$ and
taking into account (\ref{i1}) we see that
$$
g({\cal R}(s_2),s_2)=g({\cal R}(s_2),s_3)=g({\cal R}(s_3),s_3)=0.
$$
It follows that $g({\cal R}(s_i),s_j)=0$, $i,j=1,2,3$. This means
that $(M,g)$ is anti-self-dual with zero scalar curvature.

Since $g({\cal R}(U),\omega )=0$ for every vertical vector $U$, identity $(ii)$ takes the form
$$
g({\cal R}(U),X\wedge Y+JX\wedge JY)-g({\cal R}(\sigma\times U),X\wedge JY-JX\wedge Y)=0.
$$
Setting in this identity $(X,Y)=(E_1,E_3)$ and $(X,Y)=(E_3,E_1)$ we obtain
$$
g({\cal R}(U),\bar{s}_2)-g({\cal R}(\sigma\times U),\bar{s}_3)=0, \quad g({\cal R}(U),\bar{s}_2)+g({\cal
R}(\sigma\times U),\bar{s}_3)=0.
$$
This, together with (\ref{i1}), implies $g({\cal R}(U),s^-_j)=0,\, j=1,2,3$. Thus
$$
g({\cal R}(s_i),s^-_j)=0,\quad i,j=1,2,3
$$
which means that ${\cal B}=0$. We note also that, since $dim\,M=4$, $dF= \theta\wedge F$,  where
$\theta= \delta F\circ J$, so the identity $\delta F=0$ is equivalent to $dF=0$, i.e. to $(M,g,J)$
being almost K\"ahler. It follows that if $({\cal Z},h_t, {\cal J}_{\omega})\in W_1\oplus W_2\oplus
W_4$, then $(M,g,J)$ is almost K\"ahler,  anti-self-dual and Ricci flat manifold. According to
\cite[Propostion 1]{AD} these conditions are equivalent to the base manifold being K\"ahler and
Ricci flat. For such a manifold we have $\nabla F=\delta F=0$ and ${\cal }R(U)=0$ for every
vertical vector $U$, thus conditions $(i)$ -- $(iv)$ are clearly satisfied. Now the lemma follows
from Lemma \ref{K}.

\smallskip

\begin{lemma}\label{SK} $({\cal Z},h_t, {\cal J}_{\omega})\in {\cal SK}={\cal W}_1\oplus {\cal W}_2\oplus {\cal W}_3$
if and only if $(M,g,J)$ is almost K\"ahler and ${\cal
W}_{+}(\omega)=-\displaystyle{\frac{s}{6}}\omega$.
\end{lemma}
{\bf Proof}. The defining condition for the class of semi-K\"ahler manifolds is $\delta\Omega=0$.
According to Corollary~\ref{delta-Om} and (\ref {hor}), the twistor space is semi-K\"ahler if and
only if $g(\delta\omega,X)=0$, i.e. $\delta F(X)=0$ and $g({\cal R}(U),\omega\circ\pi(\sigma))=0$
for every $\sigma\in{\cal Z}$, $U\in{\cal V}_{\sigma}$, $X\in T_{\pi(\sigma)}M$. As we have
mentioned the identity $\delta F=0$ is equivalent to $dF=0$ since $dim\,M=4$. The identity $g({\cal
R}(U),\omega\circ\pi(\sigma))=0$ for all $U\in{\cal V}_{\sigma}$ holds if and only if $g({\cal
R}(\omega), s_i)=0, i=1,2,3$. This is equivalent to $\displaystyle{\frac{s}{6}}\omega+{\cal
W}_{+}(\omega)=0$.

\smallskip

\begin{lemma}\label{G1}  $({\cal Z},h_t, {\cal J}_{\omega})\in {\cal G}_1={\cal W}_1\oplus
{\cal W}_3\oplus {\cal W}_4$ if and only if the almost complex
structure $J$ is integrable.
\end{lemma}

{\bf Proof}.  $({\cal Z},h_t, {\cal J}_{\omega})$ belongs to the
class ${\cal G}_1$ when
\begin{equation}\label{G-134}
(D_A\Omega)(A,B)-(D_{{\cal J}_{\omega}A})({\cal
J}_{\omega}A,B)=0,\quad A,B\in T{\cal Z}.
\end{equation}
It follows from Proposition~\ref{D-Om} and (\ref {hor}) that this condition holds if and only if
for every $X,Y\in TM$
$$(\nabla_X F)(X,Y)-(\nabla_{JX} F)(JX,Y)=0.$$ In dimension  $4$, the
latter identity is equivalent to $J$ being integrable..

\smallskip

\begin{lemma}\label{G2} $({\cal Z},h_t, {\cal J}_{\omega})\in {\cal G}_2={\cal W}_2\oplus
{\cal W}_3\oplus {\cal W}_4$ if and only if
$$
{\cal W}_{+}(\sigma)=\displaystyle{\frac{s}{2}g(\sigma,\omega)\omega-\frac{s}{6}\sigma}
$$
for all $\sigma\in\Lambda^2_{+}TM$.
\end{lemma}
{\bf Proof}. The condition for $({\cal Z},h_t, {\cal J}_{\omega})$
to be in the class ${\cal G}_2$ is
\begin{equation}\label{234}
\displaystyle
\mathop{\mathfrak{S}}_{A,B,C}\{(D_A\Omega)(B,C)-(D_{{\cal
J}_{\omega}A}\Omega)({\cal J}_{\omega}B,C)\}=0.
\end{equation}
By Proposition~\ref{D-Om} and (\ref {hor}) this is equivalent to the
following identities

\smallskip

\noindent $(i)$ $\displaystyle
\mathop{\mathfrak{S}}_{X,Y,Z}\{(\nabla_X F)(Y,Z)-(\nabla_{JX} F)(
JY,Z)\}=0, \quad X,Y,Z\in TM$.

\medskip

\noindent $(ii)$ $g({\cal R}(U),X\wedge Y-JX\wedge JY)-g({\cal R}(\sigma\times U),X\wedge
JY+JX\wedge Y)=0$, \hspace{1in}

\smallskip

\noindent for every $\sigma\in{\cal Z}$, $U\in{\cal V}_{\sigma}$,
$X,Y\in T_{\pi(\sigma)}M$. Identity $(i)$ is always satisfied in
dimension $4$. Identity $(ii)$ gives
$$
g({\cal R}(U),s_2)-g({\cal R}(\sigma\times U),s_3)=0.
$$
Applying the latter identity for $\sigma=s_1,s_2,s_3$, it easy to see that
$$
g({\cal R}(s_i),s_j)=0~{\mbox for}~ (i,j)\neq (1,1).
$$
The curvature decomposition and the fact
that $Trace\,{\cal W}_{+}=0$ then imply
$$\displaystyle{\frac{s}{6}}+g({\cal W}_{+}(\omega),\omega)=g({\cal
R}(\omega),\omega)=\displaystyle{\frac{s}{2}}.$$ Thus the matrix of ${\cal W}_{+}$ with respect to
the basis $s_1=\omega,s_2,s_3$ is diagonal with diagonal entries
$\displaystyle{\frac{s}{3},-\frac{s}{6},-\frac{s}{6}}$.  Therefore ${\cal
W}_{+}(\sigma)=\displaystyle{\frac{s}{2}g(\sigma,\omega)\omega-\frac{s}{6}\sigma}$.

Conversely, suppose that this identity is fulfilled. Then ${\cal
R}(\sigma)=\displaystyle{\frac{s}{2}}g(\sigma,\omega)\omega+{\cal
B}(\sigma)$. It is easy to check that if
$\sigma\in\Lambda^2_{+}T_pM$ and $\tau\in\Lambda^2_{-}T_pM$, the
endomorphisms $K_{\sigma}$ and $K_{\tau}$ of $T_pM$ commute,
$K_{\sigma}\circ K_{\tau}=K_{\tau}\circ K_{\sigma}$. This implies
that, for every $X,Y\in T_pM$, the $2$-vector $X\wedge
Y-K_{\sigma}X\wedge K_{\sigma}Y$ is orthogonal to
$\Lambda^2_{-}T_pM$, so it lies in $\Lambda^2_{+}T_pM$. In
particular, $g({\cal B}(\sigma),X\wedge Y-JX\wedge JY)=0$. We also
have $g(\omega,X\wedge Y-JX\wedge JY)=0$. Thus $g({\cal
R}(\sigma),X\wedge Y-JX\wedge JY)=0$ for every $\sigma\in{\cal
Z}$, $X,Y\in T_{\pi(\sigma)}M$. It follows that condition $(ii)$
is satisfied, hence  $({\cal Z},h_t, {\cal J}_{\omega})\in {\cal
G}_2$.

\smallskip

\begin{lemma}\label{W13} $({\cal Z},h_t, {\cal J}_{\omega})\in {\cal W}_1\oplus {\cal W}_3$ if
and only if $(M,g,J)$ is K\"ahler and scalar flat.
\end{lemma}
{\bf Proof}. Note that
$${\cal W}_1\oplus {\cal W}_3= ({\cal W}_1\oplus {\cal W}_2\oplus {\cal W}_3)\cap  ({\cal W}_1\oplus {\cal W}_3\oplus
{\cal W}_4).$$ Hence it follows from Lemmas ~\ref{SK} and \ref{G1} that $({\cal Z},h_t, {\cal
J}_{\omega})\in {\cal W}_1\oplus {\cal W}_3$ if and only if $(M,g,J)$ is K\"ahler and ${\cal
W}_{+}{(\omega})=-\displaystyle{\frac{s}{6}}\omega.$ But, as we have already mentioned, it is
well-known that for K\"ahler manifolds ${\cal W}_{+}{(\omega})=\displaystyle{\frac{s}{3}}\omega$
and the above identity implies that $s=0$ . The converse statement follows from the fact that a
K\"ahler manifold is scalar flat if and only if it is anti-self-dual.

\smallskip

\begin{lemma}\label{W23} $({\cal Z},h_t, {\cal J}_{\omega})\in {\cal W}_2\oplus {\cal W}_3$ if
and only if $(M,g,J)$ is K\"ahler and scalar flat.
\end{lemma}
{\bf Proof}. It follows from Lemmas ~\ref{SK} and \ref{G2} that $({\cal Z},h_t, {\cal
J}_{\omega})\in {\cal W}_2\oplus {\cal W}_3$ if and only if $(M,g,J)$ is almost K\"ahler,
anti-self-dual and scalar flat. Now the lemma follows from Proposition 1 in \cite{AD} according to
which these conditions are equivalent to the base manifold being K\"ahler and scalar flat.

\smallskip

\begin{lemma}\label{SKcupH} $({\cal Z},h_t, {\cal J}_{\omega})\in{\cal W}_3$ if and only if
$(M,g,J)$ is a K\"ahler and scalar flat.
\end{lemma}
{\bf Proof}. The lemma follows from Lemmas ~\ref{W13} and \ref{W23}.

\smallskip

\begin{lemma}\label{H}  $({\cal Z},h_t, {\cal J}_{\omega})\in {\cal H}={\cal W}_3\oplus
{\cal W}_4$  if and only if the almost complex structure $J$ is
integrable and $${\cal
W}_{+}(\sigma)=\displaystyle{\frac{s}{2}g(\sigma,\omega)\omega-\frac{s}{6}\sigma}$$
for all $\sigma\in\Lambda^2_{+}TM$.
\end{lemma}
{\bf Proof}. The proof follows from Lemmas ~\ref{G1} and ~\ref{G2}.

\medskip

We are now ready to prove Theorem~\ref{GH}.

\smallskip

{\bf Proof of Theorem ~\ref{GH}}.

\smallskip

It follows from Lemmas \ref{K} and \ref{W124} that $${\cal K}= {\cal W}_1= {\cal W}_2= {\cal W}_4=
{\cal W}_1\oplus {\cal W}_2= {\cal W}_1\oplus {\cal W}_4={\cal W}_2\oplus {\cal W}_4={\cal
W}_1\oplus {\cal W}_2\oplus {\cal W}_4.$$

Lemmas \ref{W13}, \ref{W23} and \ref{SKcupH} imply that
$${\cal W}_3= {\cal W}_1\oplus{\cal W}_3= {\cal W}_2\oplus {\cal
W}_3.$$

Hence the first part of the theorem follows from Lemmas \ref{SK},
\ref{G1}, \ref{G2} and \ref{H}.

The statements (i)--(vi) follow respectively from Lemmas \ref{K},
\ref{SKcupH} , \ref{H}, \ref{SK}, \ref{G1} and \ref{G2} .

\smallskip

\noindent {\it Remark}. Concerning the geometric conditions in
Theorem $1$ we note that:
\smallskip

\noindent $(i).$ Any compact K\"ahler and Ricci flat surface is
either a complex torus, a hyperelliptic surface with the flat
metric, a $K3$-surface with a Calabi-Yau metric or its ${\Bbb Z}$
or ${\Bbb Z\oplus\Bbb Z}$ quotient.
\smallskip

\noindent $(ii).$ The spectrum of the anti-self-dual Weyl tensor
${\cal W}_{+}$ of a Hermitian surface $M$ is equal to
$(\displaystyle \frac{k}{3}, -\frac{k}{6},-\frac{k}{6})$, where
$k$ is the conformal scalar curvature \cite{AG}. Hence the
curvature condition in Theorem 1, (iii) implies that $k=s$, i. e.
$\delta\theta=\|\theta\|^2$, where $\theta$ is the Lee form of $M$
\cite{AG}. If $M$ is compact, then integrating this identity and
using Stock's formula, we see that $\theta = 0$, i.e. the surface
$M$ is locally conformally K\'ahler.
\smallskip

\noindent $(iii).$ We do not know non-K\'ahler examples of compact
almost K\'ahler 4-manifolds whose anti-self-dual Weyl tensor
${\cal W}_{+}$ satisfies the condition of Theorem 1, (iv).

\section{Gray-Hervella classes of the almost complex structures  ${\cal J}_\lambda^\pm$ }
As in \cite{De}, in order to define a fibre-preserving map $f:{\cal
Z}\to {\cal Z}$ in an explicit way, we shall use the stereographic
projection of every fibre ${\cal Z}_{\pi(\sigma)}$ from the point
$\omega_{\pi(\sigma)}$ onto the plane $({\Bbb
R}\omega_{\pi(\sigma)})^{\perp}$, the orthogonal complement being
taken in $\Lambda^2_{+}T_{\pi(\sigma)}M$. This stereographic
projection $\Phi_{\sigma}$ and its inverse $\Phi_{\sigma}^{-1}$ are
given by
$$
\begin{array}{c}
\Phi_{\sigma}(\tau)=\displaystyle{\frac{\tau-g(\tau,\omega_{\pi(\sigma)})\omega_{\pi(\sigma)}}{1-g(\tau,\omega_{\pi(\sigma)})}},
\quad \tau\in{\cal Z}_{\pi(\sigma)}\setminus\{\omega_{\pi(\sigma)}\}, \\[10pt]
\Phi_{\sigma}^{-1}(\zeta)=\displaystyle{\frac{2\zeta+||\zeta||^2-1}{||\zeta||^2+1}},\quad
\zeta\in ({\Bbb R}\omega_{\pi(\sigma)})^{\perp}
\end{array}
$$
The map $\Phi_{\sigma}$ is holomorphic with respect to the standard
complex structure of ${\cal Z}_{\pi(\sigma)}$ and the complex
structure of $({\Bbb R}\omega_{\pi(\sigma)})^{\perp}$ given by
$\zeta\to \omega_{\pi(\sigma)}\times\zeta$ (the latter structure is
compatible with the metric $g$ of $\Lambda^2_{+}T_{\pi(\sigma)}M$).
As usual, we also set $\Psi_{\sigma}(\sigma)=\infty$, the "ideal"
element of the plane $({\Bbb R}\omega_{\pi(\sigma)})^{\perp}$.

Let $\lambda=a+ib\in{\Bbb C}$ and set $F_\lambda(\zeta)=\lambda\zeta$ for $\zeta\in ({\Bbb
R}\omega_{\pi(\sigma)})^{\perp}$. Then  $$f_\lambda^+(\sigma)= \Phi_{\sigma}^{-1}\circ
F_\lambda\circ\Phi_{\sigma}(\sigma)$$ is a self-map of ${\cal Z}$ whose restriction to any fibre is
holomorphic. Similarly, denote by $\Psi_{\sigma}$ the stereographic projection of ${\cal
Z}_{\pi(\sigma)}$ from the point $-\omega_{\pi(\sigma)}$ onto the plane $({\Bbb
R}\omega_{\pi(\sigma)})^{\perp}$. Set $f_\lambda^-(\sigma)=\Psi_{\sigma}^{-1}\circ F_\lambda\circ
\Psi_{\sigma}(\sigma)$. In this way we obtain another  self-map  of ${\cal Z}$ whose restriction to
any fibre is anti-holomorphic. Clearly, the points $f_\lambda^-(\sigma)$ and $f_\lambda^+(\sigma)$
are symmetric with respect to the plane $({\Bbb R}\omega_{\pi(\sigma)})^{\perp}$.

The maps $f_\lambda^\pm:{\cal Z}\to {\cal Z}$ are given by the
following explicit formula:
$$
\begin{array}{c}
f_\lambda^{\pm}(\sigma)=[(a^2+b^2+1)+(a^2+b^2-1)g(\sigma,\omega_{\pi(\sigma)})]^{-1}\times\\[6pt]
\{2a\sigma-2b\sigma\times\omega_{\pi(\sigma)}-2ag(\sigma,\omega_{\pi(\sigma)})\omega_{\pi(\sigma)}\\[6pt]
\pm[(a^2+b^2-1)+(a^2+b^2+1)g(\sigma,\omega_{\pi(\sigma)})]\omega_{\pi(\sigma)}\}.
\end{array}
$$
Denote by ${\cal J}_\lambda ^\pm$  the almost complex structure on ${\cal Z}$ defined by means of
the map $f_\lambda^{\pm}(\sigma)$. Note that $f_0^\pm\equiv\mp\omega$ and ${\cal J}_0^\pm$ is the
almost complex structure on ${\cal Z}$ yielded by the almost complex structure $\mp J$ on $M$ and
discussed in the preceding section. The structure ${\cal J}_\lambda^+$ is denoted by $J_{\lambda
Id}$ in \cite{De} where the integrability condition for this structure is found when the base
manifold $(M,g,J)$ is K\"ahler. Note also that $f_1^+(\sigma)=\sigma$ and ${\cal J}_1^+$ is the
Atiyah-Hitchin-Singer almost complex structure, whereas $f_{-1}^-(\sigma)=-\sigma$ and ${\cal
J}_{-1}^-$ is the Eells-Salamon almost complex structure. The Gray-Hervella classes of these
structures have been determined in \cite{M}.

\smallskip

Since the restrictions to the fibres of the map $f_\lambda^-(\sigma)$ are not holomorphic,
Corollary~\ref{hol f} implies the following.

\begin{cor}
The almost complex structure ${\cal J}_\lambda ^-$ is never integrable.
\end{cor}




In this section we shall we discuss the possible Gray-Hervella classes of the almost Hermitian
manifolds $(Z, h_t, {\cal J}_\lambda^\pm)$. To do this we need to compute ${\cal V}(f^\pm
_\lambda)_ \ast(X^h_{\sigma})$, $X\in T_{\pi(\sigma)}M$. Taking a section $s$ of $\Lambda^2_{+}TM$
around the point $p=\pi(\sigma)$ such that $s(p)=\sigma$ and $\nabla s|_p=0$, we have ${\cal
V}(f^\pm _\lambda)_ \ast(X^h_{\sigma})=\nabla_{X}(f^\pm_\lambda\circ s)$. Using this formula, we
can get an explicit expression for ${\cal V}(f^\pm _\lambda)_ \ast(X^h_{\sigma})$ which simplifies
considerably in the case when $(M,g,J)$ is a K\"ahler manifold or when $|\lambda|=1$. In fact,
$${\cal V}(f^\pm _\lambda)_ \ast(X^h_{\sigma})=0$$
in the first case and
\begin{equation}\label{star} {\cal V}(f^\pm
_\lambda)_ \ast(X^h_{\sigma})=-b\sigma\times\nabla_{X}\omega  +(\pm
1-a)[g(\sigma,\nabla_{X}\omega)\omega_{\pi(\sigma)}
+g(\sigma,\omega)\nabla_{X}\omega]
\end{equation} in the case when $|\lambda|=1$. Let us note that if $|\lambda|=1$, say $\lambda=e^{i\theta}$, the point $f_\lambda^+(\sigma)$ is
obtained by rotating $\sigma$  around the line ${\Bbb
R}\omega_{\pi(\sigma)}$ at angle $\theta$.

We are now ready to prove the following

\smallskip

\begin{th} \label{GHC1} Let $(M,g,J)$ be a K\"ahler manifold and $\lambda\neq 0$ be a complex number.

\noindent $(i)$ The possible Gray-Hervella classes of the twistor space $({\cal Z},h_t, {\cal
J}_\lambda^-)$ are ${\cal W},~ {\cal QK}={\cal W}_1\oplus {\cal W}_2$ and $~{\cal SK}={\cal
W}_1\oplus {\cal W}_2\oplus {\cal W}_3$ .  Moreover

\smallskip

\noindent $(i_1)$ $({\cal Z},h_t, {\cal J}_\lambda^-)\in {\cal SK}={\cal W}_1\oplus {\cal
W}_2\oplus {\cal W}_3$ if and only if $(M,g,J)$ is scalar flat.

\medskip

\noindent $(i_2)$ $({\cal Z},h_t, {\cal J}_\lambda^-)\in {\cal QK}={\cal W}_1\oplus {\cal W}_2$ if
and only if $(M,g,J)$ is Ricci flat.

\medskip

\noindent $(ii)$ The possible Gray-Hervella classes of the twistor space $({\cal Z},h_t, {\cal
J}_\lambda^+)$ are ${\cal W}$ and ~${\cal W}_3={\cal SK}\cap {\cal H}$. The latter case occurs if
and only if $(M,g,J)$ is scalar flat.
\end{th}

\bigskip

{\bf Proof}. Given a point $p\in M$, we choose an orthonormal frame
of vector fields $A_1,...,A_4$ around $p$ such that $A_3=JA_2$,
$A_4=JA_1$ and use this frame to define sections $s_1,s_2,s_3$ of
$\Lambda^2_{+}TM$ via (\ref{s-basis}). Then $\omega=s_3$ and

$$
\begin{array}{c}

s_1=A_1\wedge A_2-JA_1\wedge JA_2, s_2=A_1\wedge JA_2+JA_1\wedge
A_2, s_3=A_1\wedge JA_1+A_2\wedge JA_2.

\end{array}
$$

Suppose that $(M,g,J)$ is a K\"ahler $4$-manifold. Then, as we have mentioned,
\begin{equation}\label{ki}
{\cal R}(s_1)={\cal R}(s_2)=0,\quad g({\cal
R}(s_3),s_3)=\frac{s}{3}s_3.
\end{equation}
In particular the K\"ahler metric $g$ is anti-self-dual if and only if it is scalar flat.

To determine the possible Gray-Hervella classes of the twistor space $({\cal Z},h_t, {\cal
J}_\lambda^{\pm})$ of an almost Hermitian manifold $(M,g,J)$ we shall need several technical
lemmas. Next we shall always assume that $(M,g,J)$ is a K\"ahler $4$-manifold with K\"ahler
$2$-vector $\omega$ and scalar curvature $s$ and that $\lambda\neq 0$ is an arbitrary complex
number.

\begin{lemma}\label{k-SK} $({\cal Z},h_t,{\cal J}^{\pm}_{\lambda})\in {\cal
SK}={\cal W}_1\oplus {\cal W}_2\oplus {\cal W}_3$ if and only if
$s=0$.
\end{lemma}

{\bf Proof}. We know that ${\cal V}(f^\pm _\lambda)_ \ast(X^h_{\sigma})=0$ for all $X\in
T_{\pi(\sigma)}M$, hence, by Corollary~\ref{delta-Om}, the fundamental $2$-form of $(h_t,{\cal
J}^{\pm}_{\lambda})$ is co-closed if and only if
\begin{equation}\label{k-sk-eq}
g({\cal R}(U),f^{\pm}(\sigma))=0.
\end{equation}
for every $\sigma\in{\cal Z}$ and $U\in{\cal V}_{\sigma}$. Setting
in this identity $\sigma=s_1(p)$, $U=s_3(p)$ for $p\in M$ and taking
into account (\ref{ki}), we obtain $(a^2+b^2-1)g({\cal
R}(s_3),s_3)=0$. Hence $g({\cal R}(s_3),s_3)=0$ if $a^2+b^2\neq 1$.
If $a^2+b^2=1$ we set $\sigma=\frac{1}{\sqrt 2}(s_1+s_3)$ and
$U=s_1-s_3$. Then ${\sqrt 2}f(\sigma)=as_1+bs_2 \pm s_3$ and
identity (\ref{k-sk-eq}) gives again $g({\cal R}(s_3),s_3)=0$. It
follows that $s=0$.

Conversely, if $s=0$, we have $g({\cal R}(s_i),s_j)=0$, $i,j=1,2,3$,
so $g({\cal R}(\sigma),\tau)=0$ for every
$\sigma,\tau\in\Lambda^2_{+}TM$. In particular, identity
(\ref{k-sk-eq}) is fulfilled, hence $(h_t,{\cal J}^{\pm}_{\lambda})$
is semi-K\"ahler.

\smallskip

\begin{lemma}\label{k-QK} \noindent $(i)$  $({\cal Z}, h_t,{\cal J}_\lambda^-) \in {\cal QK}={\cal W}_1\oplus {\cal W}_2$
if and only if $(M,g,J)$ is Ricci flat.

\smallskip
\noindent $(ii)$ $({\cal Z},h_t,{\cal J}_\lambda^+)$ never belongs
to the class ${\cal W}_1\oplus {\cal W}_2 $.
\end{lemma}

{\bf Proof}. Suppose that $({\cal Z}, h_t,{\cal J}_\lambda^{\pm}) \in {\cal W}_1\oplus {\cal W}_2$.
By the defining condition for the class of quisi-K\"ahler manifolds
\begin{equation}\label{q0}
(D_{A}\Omega)(B,C)+(D_{{\cal J}_\lambda^{\pm}A}\Omega)({\cal J}_\lambda^{\pm}B,C)=0,\quad A,B\in
T{\cal Z}.
\end{equation}
Hence, according to the second formula of Proposition~\ref{D-Om},
\bigskip
\begin{equation}\label{q1}
\begin{array}{c}
g({\cal R}(U),X\wedge Y+ K_{f_\lambda ^{\pm}(\sigma)}X\wedge K_{f_\lambda^{\pm}(\sigma)}Y)\\[6pt]
-g({\cal R}(\sigma\times U),X\wedge K_{f_\lambda^{\pm}(\sigma)}Y-K_{f_\lambda^{\pm}(\sigma)}X\wedge
Y)=0.
\end{array}
\end{equation}
for every $\sigma\in{\cal Z}$, $X,Y\in T_{\pi(\sigma)}M$, $U\in {\cal V}_{\sigma}$. Setting
$Y=K_{f_\lambda^{\pm}(\sigma)}X$ we get

\begin{equation}\label{q2}
g({\cal R}(U),X\wedge K_{f_\lambda^{\pm}(\sigma)}X)=0,
\end{equation}
or, equivalently,
\begin{equation}\label{q3}
g({\cal R}(U),X\wedge K_{f_\lambda^{\pm}(\sigma)}Y-K_{f_\lambda^{\pm}(\sigma)}X\wedge Y)=0.
\end{equation}
It is easy to check by means of (\ref{com}) that for any $\tau\in \Lambda^2_{+}T_pM$ and $X,Y\in
T_pM$ with $X\perp Y$, the $2$-vector $X\wedge K_{\tau}Y-K_{\tau}X\wedge Y$ is orthogonal to
$\Lambda^2_{+}T_pM$, hence it lies in $\Lambda^2_{-}T_pM$. Moreover, for every $\tau\in
\Lambda^2_{+}T_pM$, every vector of $\Lambda^2_{-}T_pM$ is a linear combination of vectors of the
form $X\wedge K_{\tau}Y-K_{\tau}X\wedge Y$ with $X\perp Y$ and vectors of the form $Z\wedge
K_{\tau}Z$, \, $X,Y,Z\in T_pM$. Indeed, if $a_1,...,a_4$ is an orthonormal basis of $T_pM$ such
that $a_3=K_{\tau}a_2$, $a_4=K_{\tau}a_1$, then it is positively oriented and $a_1\wedge
a_2-a_3\wedge a_4=-(a_1\wedge K_{\tau}a_3-K_{\tau}a_1\wedge a_3)$, $a_1\wedge a_3-a_4\wedge
a_2=a_1\wedge K_{\tau}a_2-K_{\tau}a_1\wedge a_2$, $a_1\wedge a_4-a_2\wedge a_3=a_1\wedge
K_{\tau}a_1-a_2\wedge K_{\tau}a_2$. Thus it follows from (\ref{q2}) and (\ref{q3}) that
$$
g({\cal R}(U),s^{-})=0
$$
for every $\sigma\in{\cal Z}$, $U\in{\cal V}_{\sigma}$, $s^{-}\in \Lambda^2_{-}T_{\pi(\sigma)}M$.
In particular, $g({\cal R}(s_3),s^{-})=0$, hence in view of (\ref{ki}), $ {\cal
R}(s_3)=\displaystyle{\frac{s}{3}}s_3$.  Now, setting $\sigma=s_1(p)$ and $U=s_3(p)$, $p\in M$, in
$(\ref{q3})$, we obtain
\begin{equation}\label{s3}
sg(s_3,X\wedge Y+ K_{f_\lambda^{\pm}(s_1)}X\wedge K_{f_\lambda ^{\pm}(s_1)}Y)=0.
\end{equation}
This identity for $(X,Y)=(A_1,A_2)$ and $(X,Y)=(A_1,A_3)$ gives
$$
sa(a^2+b^2-1)=0,\quad sb(a^2+b^2-1)=0.
$$
Hence $s=0$ if $a^2+b^2\neq 1$. If $a^2+b^2=1$, we set $\sigma =\frac{1}{\sqrt 2}(s_1+s_3)$ and
$U=s_1-s_3$. We have ${\sqrt 2}f(\sigma)=as_1+bs_2 \pm s_3$ and identity (\ref{s3}) with
$(X,Y)=(A_1,A_2)$ and $(A_1,A_3)$ gives $as=0$ and $bs=0$. Therefore $s=0$. It follows that ${\cal
R}(\tau)=0$ for every $\tau\in\Lambda^2_{+}TM$. Since $(M,g,J)$ is a K\"ahler manifold this is
equivalent to the metric $g$ being Ricci flat. Moreover, in view of the third formula of
Proposition~\ref{D-Om}, we get from (\ref{q0}) that
\begin{equation}\label{q-eps}
g((f_\lambda^{\pm})_{\ast}(U),Y\wedge
Z)+g((f_\lambda^{\pm})_{\ast}({\cal
J}_\lambda^{\pm}U),K_{f_\lambda^{\pm}(\sigma)}Y\wedge Z)=0
\end{equation}
for $Y,Z\in T_{\pi(\sigma)}M$ and $U\in{\cal V}_{\sigma}$. The
restriction of $f_\lambda^+$ to any fibre of ${\cal Z}$ is
holomorphic, hence
$$
(f_\lambda^{+})_{\ast}({\cal J}_\lambda^{+}U)={\cal J}_\lambda^{+}(f_\lambda
^{+})_{\ast}(U)=f_\lambda^{+}(\sigma)\times (f_\lambda^{+})_{\ast}(U).
$$
Then, by (\ref{vp})
$$
g((f_\lambda^{+})_{\ast}({\cal
J}_\lambda^{+}U),K_{f_\lambda^{+}(\sigma)}Y\wedge Z)=g((f_\lambda
^{+})_{\ast}(U),Y\wedge Z).
$$
This and (\ref{q-eps}) imply $ (f_\lambda^{+})_{\ast}(U)=0$.
Therefore the restriction of $f_\lambda^{+}$ to every fibre of
${\cal Z}$ is a constant map which is a contradiction.
\smallskip

Now suppose that the K\"ahler manifold $M$ is Ricci flat. In this case ${\cal R}(\tau)=0$ for every
$\tau\in\Lambda^2_{+}TM$. Then, in view of Proposition~\ref{D-Om}, in order to prove that $({\cal
Z}, h_t,{\cal J}_\lambda^-) \in {\cal W}_1\oplus {\cal W}_2$ it is enough to show that, for every
$U\in{\cal V}_{\sigma}$ and $Y,Z\in T_{\pi(\sigma)}M$, we have
$$
(D_{U}\Omega)(Y^h_{\sigma},Z^h_{\sigma})+(D_{{\cal
J}_\lambda^{-}U}\Omega)((K_{f_\lambda^{-}(\sigma)}Y)^h_{\sigma},Z^h_{\sigma})=0
$$
This is equivalent to identity (\ref{q-eps}) for $f_\lambda^{-}$.
The restriction of the map $f_\lambda^{-}$  to any fibre of ${\cal
Z}$ is anti-holomorphic, hence
$$
(f_\lambda^-)_{\ast}({\cal J}_\lambda^-U)=-{\cal
J}_\lambda^-(f_\lambda^-)_{\ast}(U)=-f_\lambda^-(\sigma)\times (f_\lambda^-)_{\ast}(U).
$$
This, in view of (\ref{vp}), implies that identity (\ref{q-eps}) for
$f_\lambda^{-}$ is fulfilled. Therefore $({\cal Z}, h_t,{\cal
J}_\lambda^-) \in {\cal W}_1\oplus {\cal W}_2$.

\smallskip

\begin{lemma}\label{sum3} \noindent $(i)$ $({\cal Z},h_t,{\cal
J}^{-}_{\lambda})\in {\cal W}_1\oplus {\cal W}_2\oplus {\cal W}_4$
if and only if $(M,g,J)$ is Ricci flat.

\smallskip

\noindent $(ii)$ $({\cal Z},h_t,{\cal J}^{+}_{\lambda})$ never
belongs to the class ${\cal W}_1\oplus {\cal W}_2\oplus {\cal W}_4$.
\end{lemma}

{\bf Proof}.  Suppose that $({\cal Z},h_t,{\cal J}^{\pm}_{\lambda})$
is of class ${\cal W}_1\oplus {\cal W}_2\oplus {\cal W}_4$. Then
$$
(D_{X^h_{\sigma}}\Omega)(X^h_{\sigma},U)+(D_{{\cal
J}^{\pm}_{\lambda}X^h_{\sigma}}\Omega)({\cal
J}^{\pm}_{\lambda}X^h_{\sigma},U)=
-\frac{1}{2}||X||^2\delta\Omega(U)
$$
for $X\in T_{\pi(\sigma)}M$, $U\in{\cal V}_{\sigma}$.  By
Proposition~\ref{D-Om} and Corollary~\ref{delta-Om}, this is
equivalent to
\begin{equation}\label{k-124-eq1}
2g({\cal R}(U),X\wedge K_{f_\lambda^{\pm}(\sigma)}X)=||X||^2g({\cal
R}(U),f_\lambda^{\pm}(\sigma))
\end{equation}
or
\begin{equation}\label{k-124-eq2}
g({\cal R}(U),X\wedge
K_{f_\lambda^{\pm}(\sigma)}Y-K_{f_\lambda^{\pm}(\sigma)}X\wedge
Y)=g(X,Y)g({\cal R}(U),f_\lambda^{\pm}(\sigma)).
\end{equation}
for $X,Y\in T_{\pi(\sigma)}M$, $U\in{\cal V}_{\sigma}$. Take an
orthonormal  basis $E_1,...,E_4$ of $T_{\pi(\sigma)}M$ such
$E_3=K_{f_\lambda ^{\pm}(\sigma)}E_2$,
$E_4=K_{f_\lambda^{\pm}(\sigma)}E_1$. Then identity
(\ref{k-124-eq1}) for $X=E_1$ and $X=E_2$ gives
$$
2g({\cal R}(U),E_1\wedge E_4)=g({\cal R}(U),f_\lambda^{\pm}(\sigma)),\quad 2g({\cal R}(U),E_2\wedge
E_3)=g({\cal R}(U),f_\lambda^{\pm}(\sigma)).
$$
These identities imply $g({\cal R}(U), E_1\wedge E_4-E_2\wedge E_3)=0$.
 Moreover, setting in (\ref{k-124-eq2}) $(X,Y)=(E_1,E_3)$ and
$(X,Y)=(E_1,E_2)$ we get
$$
g({\cal R}(U),E_1\wedge E_2-E_3\wedge E_4)=0,\quad g({\cal
R}(U),E_1\wedge E_3-E_4\wedge E_2)=0.
$$
It follows that $g({\cal R}(U),s^{-})=0$ for every $s^{-}\in
\Lambda^2_{-}T_{\pi(\sigma)}M$, hence $g({\cal R}(\sigma),s^{-})=0$
for $\sigma\in{\cal Z}$ and $s^{-}\in
\Lambda^2_{-}T_{\pi(\sigma)}M$.  This and (\ref{ki}) imply
$$
{\cal R}(\sigma)=\frac{s}{3}g(\sigma,\omega)\omega, \quad
\sigma\in\Lambda^2_{+}TM.
$$
Then, by Corollary~\ref{delta-Om},
$$
\delta\Omega(s_3)=-tg({\cal R}(s_1\times
s_3),f_\lambda^{\pm}(s_1))=t\frac{s}{3}g(s_2,s_3)g(s_3,f_\lambda^{\pm}(s_1)=0.
$$
Moreover, $\delta\Omega(s_1)=tg({\cal R}(s_2),f_\lambda^{\pm}(s_3))=0$ and
$\delta\Omega(s_2)=-tg({\cal R}(s_1),f_\lambda^{\pm}(s_3))=0$. It follows that $\delta\Omega=0$,
hence $s=0$ by Lemma~\ref{k-SK}. Finally note that an almost Hermitian manifold with
$\delta\Omega=0$ belongs to the class ${\cal W}_1\oplus {\cal W}_2\oplus {\cal W}_4$ if and only if
it belongs to the class ${\cal W}_1\oplus {\cal W}_2$ . Hence the lemma follows from
Lemmas~\ref{k-QK} and \ref{k-SK}.

\begin{lemma}\label{Her 1} \rm ({\cite{De})}   $({\cal Z}, h_t,{\cal
J}_\lambda^+) \in {\cal H}={\cal W}_3\oplus {\cal W}_4$ if and only
if $(M,g,J)$ is scalar flat.
\end{lemma}

{\bf Proof}. By Corollaries~\ref{hol f} and ~\ref{N}, the almost
complex structure ${\cal J}_\lambda^+$ is integrable if and only if
\begin{equation}\label{int0}
\begin{array}{c}
g({\cal R}(X\wedge K_{f_\lambda^+(\sigma)}Y+K_{f_\lambda^+(\sigma)}X\wedge Y),U)\\[6pt]
+g({\cal R}(X\wedge Y-K_{f_\lambda^+(\sigma)}X\wedge K_{f_\lambda^+(\sigma)}Y) ,\sigma\times U)=0
\end{array}
\end{equation}
It is easy to check that for every $\tau\in\Lambda^2_{+}T_pM$ and
$X,Y\in T_pM$, the $2$-vector $X\wedge K_{\tau}Y+K_{\tau}X\wedge
Y\in \Lambda^2_{+}T_pM$ (and is orthogonal to $\tau$). Therefore, in
view of (\ref{ki}), ${\cal J}_\lambda^+$ is integrable if and only
if
\begin{equation}\label{int}
\begin{array}{c}
g(X\wedge K_{f_\lambda^+(\sigma)}Y+K_{f_\lambda^{+}(\sigma)}X\wedge Y,s_3)g({\cal R}(s_3),U)\\[6pt]
+g(X\wedge Y-K_{f_\lambda^+(\sigma)}X\wedge K_{f_\lambda^+(\sigma)}Y,s_3)g({\cal
R}(s_3),\sigma\times U)=0
\end{array}
\end{equation}
for $X,Y\in T_{\pi(\sigma)}M$ and $U\in{\cal V}_{\sigma}$. Set $\sigma=s_1$ and  $U=s_3$. Then,
since ${\cal R}(s_2)=0$, identity (\ref{int}) becomes
\begin{equation}\label{int-1}
g(X\wedge K_{f_\lambda^+(\sigma)}Y+K_{f_\lambda^{+}(\sigma)}X\wedge Y,s_3)g({\cal R}(s_3),s_3).
\end{equation}
For $(X,Y)=(A_1,A_2)$ and $(X,Y)=(A_1,A_3)$, the vector $X\wedge
K_{f_\lambda^+(\sigma)}Y+K_{f_\lambda^+(\sigma)}X\wedge Y$ is collinear to $-2bs_3+(a^2+b^2-1)s_2$
and $2as_3-(a^2+b^2-1)s_1$, respectively. Then identity (\ref{int-1}) gives
$$
bg({\cal R}(s_3),s_3)=0, \quad  ag({\cal R}(s_3),s_3)=0.
$$
Therefore $g({\cal R}(s_3),s_3)=0$, thus $s=0$.  This shows that if
${\cal J}_\lambda ^+$ is integrable, then $(M,g,J)$ is  scalar flat.

Conversely, suppose that $(M,g,J)$ is K\"ahler and scalar flat. Then ${\cal
V}(f_\lambda^+)_{\ast}(X^h_{\sigma})=0$ for every $\sigma\in{\cal Z}$ and $X\in T_{\pi(\sigma)}M$.
Hence, by Corollary~\ref{N}, ${\cal H}N(X^h,Y^h)=0$ for every $X,Y$.  We also have $g({\cal
R}(s_i),s_j)=0$, $i,j=1,2,3$ since $s=0$. Thus $g({\cal R}(\sigma),\tau)=0$ for every
$\sigma,\tau\in\Lambda^2_{+}TM$. Recall that for every $\tau\in\Lambda^2_{+}TM$ and $X,Y\in
T_{\pi(\tau)}M$, the $2$-vector $X\wedge K_{\tau}Y+K_{\tau}X\wedge Y$ lies in $\Lambda^2_{+}TM$.
Then by Corollary~\ref{N} we get ${\cal V}N(X^h,Y^h)=0$. Finally, the map $f_\lambda^+$ is
holomorphic, hence, by Corollary~\ref{hol f} we have ${\cal H}(N(X^h_{\sigma},U))=0$ for every
$U\in{\cal V}_{\sigma}$. Now Corollary~\ref{N} implies that $N=0$.

\smallskip

\begin{lemma}\label{k-G1} \noindent $(i)$  $({\cal Z},h_t,{\cal J}^{-}_{\lambda})$ never belongs to the class
 ${\cal G}_1={\cal W}_1\oplus {\cal W}_3\oplus {\cal W}_4$.

\smallskip

\noindent $(ii)$ $({\cal Z}, h_t,{\cal J}^{+}_{\lambda})\in {\cal
G}_1={\cal W}_1\oplus {\cal W}_3\oplus {\cal W}_4$ if and only if
$(M,g,J)$ is scalar flat.

\end{lemma}

{\bf Proof}. By the definition of the class ${\cal G}_1$ (\cite{GH}),  $({\cal Z},h_t,{\cal
J}^{\pm}_{\lambda})\in {\cal G}_1$ if and only if

$$h_t(N(A,B),C)+h_t(N(C,B),A)=0$$
for all $A,B,C\in T{\cal Z}$. By Corollary \ref{N} this is
equivalent to the identity

\begin{equation}\label{k-g1}
\begin{array}{c}
tg({\cal R}(U),X\wedge K_{f_\lambda^{\pm}(\sigma)}Y+ K_{f_\lambda^{\pm}(\sigma)}X\wedge Y)\\[6pt]
+ tg({\cal R}(\sigma\times U),X\wedge Y-K_{f_\lambda^{\pm}(\sigma)}X\wedge
K_{f_\lambda^{\pm}(\sigma)}Y)\\[6pt]
+2g((f_\lambda^{\pm})_{\ast}(\sigma\times U)-f_\lambda^{\pm} (\sigma)\times
(f_\lambda^{\pm})_{\ast}(U),X\wedge Y)=0.
\end{array}
\end{equation}

To prove $(i)$ note that the restriction of $f_\lambda^- $ on the fibre is anti-holomorphic, thus
$f_\lambda^{\pm} (\sigma)\times (f_\lambda^{\pm})_{\ast}(U)=- (f_\lambda^{\pm})_{\ast}(\sigma\times
U)$. Hence, if $({\cal Z}, h_t,{\cal J}^{-}_{\lambda})\in {\cal G}_1$, then, setting
$\sigma=s_3(p)$, $U=s_1(p)$, and taking into account that ${\cal R}(s_1)={\cal R}(s_2)=0$, we
obtain from (\ref{k-g1}) $(f_\lambda ^{-})_{\ast\,s_3}(s_1)=0$. But a straightforward computation
shows that $(f_\lambda^{-})_{\ast\,s_3}(s_1)\neq 0$, a contradiction.

To prove $(ii)$ notice that the restriction of $f_\lambda^+ $ on the fibre is holomorphic and
identity (\ref{k-g1}) takes the form  (\ref{int0}). Hence $({\cal Z},h_t,{\cal J}^{+}_{\lambda})$
is of class ${\cal G}_1$ if and only if its of class ${\cal H}$, the later condition being
equivalent to $s=0$ by Lemma~\ref{Her 1}.

\smallskip

\begin{lemma}\label{k-G2}

\noindent $(i)$  $({\cal Z}, h_t,{\cal J}^{-}_{\lambda})$ never belongs to the class ${\cal
G}_2={\cal W}_2\oplus {\cal W}_3\oplus {\cal W}_4$.

\smallskip

\noindent $(ii)$ $({\cal Z}, h_t,{\cal J}^{+}_{\lambda})\in {\cal
G}_2={\cal W}_2\oplus {\cal W}_3\oplus {\cal W}_4$  if an only if
$(M,g,J)$ is scalar flat.
\end{lemma}

{\bf Proof}. The structure $(h_t,{\cal J}^{\pm}_{\lambda})$ is of
class ${\cal G}_2$ if and only if  \cite{GH}
$$
\displaystyle \mathop{\mathfrak{S}}_{A,B,C} h_t(N(A,B),{\cal J}^{\pm}_{\lambda}C)=0,\quad A,B,C\in
T{\cal Z}. $$ For $A=X^h_{\sigma}$, $B=U\in {\cal V}_{\sigma}$, $C=Y^h_{\sigma}$ this identity and
Corollary~\ref{N} imply
\begin{equation}\label{kg2}
\begin{array}{c}
2g({\cal J}^{\pm}_{\lambda}(f_\lambda^{\pm})_{\ast}(U)- (f_\lambda^{\pm})_{\ast}({\cal J}^{\pm}_{\lambda}U),
X\wedge K_{f_\lambda^{\pm}(\sigma)}Y+K_{f_\lambda^{\pm}(\sigma)}X\wedge Y)\\[6pt]
-tg({\cal R}(X\wedge K_{f_\lambda^{\pm}(\sigma)}Y+K_{f_\lambda^{\pm}(\sigma)}X\wedge Y),\sigma\times U)\\[6pt]
+tg({\cal R}(X\wedge Y-K_{f_\lambda^{\pm}(\sigma)}X\wedge K_{f_\lambda^{\pm}(\sigma)}Y,U)=0.
\end{array}
\end{equation}
Set $\sigma=\omega(p)$, $p\in M$. We have $f_\lambda^{\pm}(\omega)=\pm\omega$, so
$K_{f_\lambda^{\pm}(\omega)}=\pm J$.  Then, setting $(X,Y)=(A_1,A_3)$, $(X,Y)=(A_1,A_2)$ and taking
into account (\ref{ki}), we obtain form (\ref{kg2}) that
$$
g({\cal J}^{\pm}_{\lambda}(f_\lambda^{\pm})_{\ast}(U)- (f_\lambda^{\pm})_{\ast}({\cal
J}^{\pm}_{\lambda}U),s_i)=0, \quad i=1,2 .$$ The map $f_\lambda^-$ is anti-holomorphic on the
fibres of ${\cal Z}$, so the latter identity gives
$$
g((f_\lambda^-)_{\ast,\,s_3}({\cal J}^{\pm}_{\lambda}U),s_i)=0,\quad i=1,2, \quad U\in{\cal
V}_{s_3}.
$$ We set $U=s_2(p)$
and $U=s_1(p)$ and compute $$
f^{-}_{\ast,\,s_3}(s_1)=\frac{2as_1+2bs_2-(a^2+b^2-1)s_3}{2(a^2+b^2)},
f^{-}_{\ast,\,s_3}(s_2)=\frac{2as_2-2bs_1-(a^2+b^2-1)s_3}{2(a^2+b^2)}.
$$
It follows that $a=0$, $b=0$, which contradicts to the assumption
$\lambda\neq 0$. This proves statement $(i)$.

Now suppose that $({\cal Z},h_t,{\cal J}^{+}_{\lambda})$ is of class
${\cal G}_2$. Then identity (\ref{kg2}) becomes
\begin{equation}\label{kg2-1}
\begin{array}{c}
g({\cal R}(X\wedge K_{f_\lambda^{+}(\sigma)}Y+K_{f_\lambda^{+}(\sigma)}X\wedge Y),\sigma\times U)\\[6pt]
-g({\cal R}(X\wedge Y-K_{f_\lambda^{+}(\sigma)}X\wedge K_{f_\lambda^{+}(\sigma)}Y,U)=0.
\end{array}
\end{equation}
We have $f_{\lambda}^{+}(s_1)=(a^2+b^2+1)^{-1}(2as_1+2bs_2+cs_3)$ where $c=a^2+b^2-1$. Then
$$
A_1\wedge K_{f_\lambda^{+}(s_1)}A_2 + K_{f_\lambda^{+}(s_1)}A_1\wedge
A_2=(a^2+b^2+1)^{-1}(-2bs_3+cs_2).
$$
Thus, setting $\sigma=s_1$, $(X,Y)=(A_1,A_2)$, $U=s_2$ in (\ref{kg2-1}) and taking into account
that ${\cal R}(s_1)={\cal R}(s_2)=0$, we obtain $bg({\cal R}(s_3),s_3)=0$. Similarly, since
$f_{\lambda}^{+}(s_2)=(a^2+b^2+1)^{-1}(-2bs_1+2as_2+cs_3)$, setting $\sigma=s_2$,
$(X,Y)=(A_1,A_2)$, $U=s_1$ we get $ag({\cal R}(s_3),s_3)=0$. It follows $g({\cal R}(s_3),s_3)=0$,
hence $s=0$. By Lemma~\ref{Her 1} $s=0$ if and only if the almost complex structure ${\cal
J}^{+}_{\lambda}$ is integrable. In particular, if $s=0$,  $({\cal Z},h_t,{\cal J}^{+}_{\lambda})$
is of class ${\cal G}_2$ and $(ii)$ is proved.

We are now ready to prove Theorem ~\ref{GHC1}.

\smallskip

{\bf Proof of Theorem ~\ref{GHC1}}.

\smallskip

$(i)$. It follows from statements $(i)$ of Lemmas \ref{k-G1} and \ref{k-G2}, and \cite[Table I]{GH}
that the possible nontrivial Gray-Hervella classes of $({\cal Z}, h_t,{\cal J}^{-}_{\lambda})$ are
subclasses of ${\cal W}_1\oplus {\cal W}_2\oplus {\cal W}_3$ or ${\cal W}_1\oplus {\cal W}_2\oplus
{\cal W}_4.$ Moreover statements $(i)$ of Lemmas \ref{k-QK} and \ref{sum3} imply that
$${\cal W}_1\oplus{\cal W}_2= {\cal W}_1\oplus{\cal W}_2\oplus {\cal
W}_3.$$

Hence the first part of the theorem follows from statements $(i)$
of Lemmas \ref{k-SK}, \ref{k-G1} and \ref{k-G2}.

$(ii)$. Using statements $(ii)$ of Lemmas \ref{k-SK}-\ref{k-G2} we prove the second part of the
theorem in a similar way.

\medskip
Now we shall discuss the case when $|\lambda|=1$. In this case we have the simple formula
(\ref{star}) for ${\cal V}(f^\pm _\lambda)_ \ast(X^h_{\sigma})$, $\sigma\in{\cal Z}$, $X\in
T_{\pi(\sigma)}M$. This simplifies the computations that should be done in order to determine the
possible Gray-Hervella types of the almost Hermiian manifold $(h_t,{\cal J}^{\pm}_{\lambda})$. Here
we shall address only a few of the basic classes.

\begin{prop}\label{ls}
Suppose that $|\lambda|=1$ and $Re(\lambda)\neq 0,\pm 1$ for ${\cal J}^{\pm}_{\lambda}$. Then:

\smallskip
\noindent $(i)$ The almost Hermitian structure $(h_t,{\cal
J}^{-}_{\lambda})$ on the twistor space ${\cal Z}$ is
(non-integrable) quasi-K\"ahler if and only if $M$ is Ricci flat.

\smallskip
\noindent $(ii)$ The structure $(h_t,{\cal J}^{+}_{\lambda})$ is never quasi K\"ahler.

 \smallskip

\noindent  $(iii)$  The structures $(h_t,{\cal J}^{\pm}_{\lambda})$
are not nearly K\"ahler or almost K\"ahler.

\end{prop}

{\bf Proof}.  It is convenient to prove first the following.

\begin{lemma}\label{l}
 If $(h_t,{\cal J}^{\pm}_{\lambda})$ is quasi K\"ahler, then $(M,g,J)$ is K\"ahler and Ricci flat.
\end{lemma}

{\bf Proof of the lemma}.  Let $p\in M$, $X,Y,Z\in T_pM$. We have $f_{\lambda}
^{\pm}(\omega(p))=\pm (\omega(p))$, hence by Proposition~\ref{D-Om} and (\ref{star})
$$
\begin{array}{c}
0=\frac{1}{2}[(D_{X^h_{\omega(p)}}\Omega)(Y^h_{\omega(p)},Z^h_{\omega(p)})+
(D_{(JX)^h_{\omega(p)}}\Omega)((JY)^h_{\omega(p)},Z^h_{\omega(p)})]=\\[6pt]
-bg(\omega\times\nabla_{X}\omega,Y\wedge Z)+(\pm 1-a)g(\nabla_{X}\omega,Y\wedge Z)\\[6pt]
  -bg(\omega\times\nabla_{JX}\omega,JY\wedge Z)+(\pm 1-a)g(\nabla_{JX}\omega,JY\wedge Z)
\end{array}
$$
This and identity (\ref{com}) give
$$
\begin{array}{c}
-bg(\nabla_{X}J)(JY),Z)+(\pm 1-a)g((\nabla_{X}J)(Y),Z)\\[6pt]
+bg(\nabla_{JX}J)(Y),Z)+(\pm 1-a)g((\nabla_{JX}J)(JY),Z).
\end{array}
$$
Thus
$$
bJ[(\nabla_{X}J)(Y)-J(\nabla_{JX}J)(Y)]+(\pm
1-a)[(\nabla_{X}J)(Y)-J(\nabla_{JX}J)(Y)]=0.
$$
By assumption $a\neq 1$ when considering ${J_{+}}$ and $a\neq -1$
for  ${\cal J}_{-}$. Thus $b^2+(\pm 1-a)^2\neq 0$ and the latter
equation implies
$$
(\nabla_{X}J)(Y)-J(\nabla_{JX}J)(Y)=0. $$ This means that the almost Hermitian structure $(g,J)$ on
$M$ is quasi K\"ahler. It follows that it is K\"ahler since $dim\,M=4$. The assumption that ${\cal
J}_{\pm}$ is quasi K\"ahler implies identity (\ref{q1}) and,  as in the proof of Lemma~\ref{k-QK},
we see that $s=0$.  Thus $M$ is K\"ahler and Ricci flat.

\smallskip

Now we are ready to prove Proposition~\ref{ls}.

\smallskip

\noindent $(i)$. If $(h_t,{\cal J}^{-}_{\lambda})$ is quasi K\"ahler, $M$ is K\"ahler and Ricci
flat by the lemma. Conversely, if $M$ is such a manifold,  $(h_t,{\cal J}^{-}_{\lambda})$ is quasi
K\"ahler by Lemma~\ref{k-QK}.

\smallskip

\noindent $(ii)$.  This statement follows form the lemma and Lemma~\ref{k-QK}.

\smallskip

\noindent $(iii)$. If $(h_t,{\cal J}^{\pm}_{\lambda})$ is nearly K\"ahler or almost K\"ahler, it is
quasi K\"ahler, hence $M$ is K\"ahler by the lemma. Then, according to Lemma~\ref{sum3} $(ii)$,
$(h_t,{\cal J}^{+}_{\lambda})$ does not belong to the class ${\cal NK}={\cal W}_1$ or to the class
${\cal AK}={\cal W}_2$. Also $(h_t,{\cal J}^{-}_{\lambda})$ is not of class ${\cal NK}$  or ${\cal
AK}$ by Lemmas~\ref{k-G1} $(i)$ and \ref{k-G2} $(i)$.

\begin{prop}
Let $|\lambda|=1$ and $Re(\lambda)\neq 0,1$.  Then the almost complex structure ${\cal
J}^{+}_{\lambda}$ is integrable if and only if $(M,g,J)$ is K\"ahler and scalar flat. In this case
$(h_t,{\cal J}^{+}_{\lambda})\in {\cal W}_3={\cal SK}\cap {\cal H}$.
\end{prop}
{\bf Proof}. Suppose that the almost complex structure ${\cal J}^{+}_{\lambda}$ is integrable. Let
$p\in M$. For $\sigma=s_1(p)$ and $X\in T_pM$, we have
$$
\begin{array}{c}
\sigma\times\nabla_{X}\omega=g(s_2,s_1\times\nabla_{X}s_3)s_2+g(s_3,s_1\times\nabla_{X}s_3)s_3=\\[6pt]
-g(s_3,\nabla_{X}s_3)s_2+g(s_2,\nabla_{X}s_3)s_3=g(s_2,\nabla_{X}s_3)s_3.
\end{array}
$$
Thus by (\ref{star})
$$
{\cal V}(f^{+} _\lambda)_ \ast(X^h_{s_1(p)})=[-bg(\nabla_{X}s_3,s_2)+(1-a)g(\nabla_{X}s_3,s_1)]s_3.
$$
It is convenient to set
$$
\phi_i=-bg(\nabla_{A_i}s_3,s_2)+(1-a)g(\nabla_{A_i}s_3,s_1)_p
$$
We have $K_{f^{+}_{\lambda}(s_1)}=aK_{s_1}+bK_{s_2}$ and, using Corollary~\ref{N}, it is easy to
see that
$$
h_t(N(A_1^h,A_2^h),A_1^h)_{s_1(p)}=2b\phi_1,\quad h_t(N(A_1^h,A_3^h),A_1^h)_{s_1(p)}=-2a\phi_1.
$$
It follows that $\phi_1=0$ since $N=0$ and $a^2+b^2\neq 0$. We also
have
$$
\begin{array}{l}
h_t(N(A_1^h,A_2^h),A_2^h)_{s_1(p)}=2b\phi_2,\quad h_t(N(A_2^h,A_4^h),A_2^h)_{s_1(p)}=2a\phi_2,\\[6pt]
h_t(N(A_1^h,A_3^h),A_3^h)_{s_1(p)}=2a\phi_3,\quad h_t(N(A_3^h,A_4^h),A_3^h)_{s_1(p)}=2b\phi_3\\[6pt]
h_t(N(A_2^h,A_4^h),E_2^h)_{s_1(p)}=2a\phi_4,\quad h_t(N(A_3^h,A_4^h),A_4^h)_{s_1(p)}=2b\phi_4.
\end{array}
$$
It follows that $\phi_2=\phi_3=\phi_4=0$. Thus
\begin{equation}\label{psi1}
-bg(\nabla_{A_i}s_3,s_2)+(1-a)g(\nabla_{A_i}s_3,s_1)=0,\quad
i=1,...,4.
\end{equation}
Now set $\sigma=s_2(p)$. We have
$$
{\cal
V}f_{\ast}(X^h_{s_2(p)})=bg(\nabla_{X}s_3,s_1)+(1-a)g(\nabla_{X}s_3,s_2)s_3
$$
and $K_{f^{+}_{\lambda}(s_2)}=-bK_{s_1}+aK_{s_2}$. Then a similar computation as above gives
\begin{equation}\label{psi2}
bg(\nabla_{A_i}s_3,s_1)+(1-a)g(\nabla_{A_i}s_3,s_2)=0,\quad
i=1,...,4.
\end{equation}
It follows from (\ref{psi1}) and (\ref{psi2}) that
\begin{equation}\label{k}
g(\nabla_{A_i}s_3,s_1)=g(\nabla_{A_i}s_3,s_2)=0,\quad i=1,...,4,
\end{equation}
Also we have $g(\nabla_{A_i}s_3,s_3)=0$ since $s_3$ is of constant length, hence, by (\ref{k}), the
almost complex structure $J$ is K\"ahlerian. By Corollary~\ref{N}, the identity ${\cal
V}N(X^h_{\sigma},Y^h_{\sigma})=0$ is equivalent to
$$
\begin{array}{c}
g({\cal R}(X\wedge K_{f^{+}_{\lambda}(\sigma)}Y+K_{f^{+}_{\lambda}(\sigma)}X\wedge Y),U)\\[6pt]
+g({\cal R}(X\wedge Y-K_{f^{+}_{\lambda}(\sigma)}X\wedge K_{f^{+}_{\lambda}(\sigma)}Y),\sigma\times
U)=0
\end{array}
$$
for $X,Y\in T_{\pi(\sigma)}M$ and $U\in{\cal V}_{\sigma}$. Setting in the latter identity
$\sigma=s_1(p)$,  $U=s_3(p)$, $(X,Y)=(A_1,A_2)$ and $(X,Y)=(A_1,A_3)$, and taking into account
(\ref{ki}), we get
$$
bg({\cal R}(s_3),s_3)=0, \quad  ag({\cal R}(s_3),s_3)=0.
$$
Therefore $g({\cal R}(s_3),s_3)=0$, thus $s=0$. This proves that if ${\cal J}^{+}_{\lambda}$ is
integrable, then $(M,g,J)$ is K\"ahler and scalar flat. The converse follows from Theorem
~\ref{GHC1} $(ii)$.

\end{document}